\documentclass[11pt]{article}

\newcommand{\real}{{\bf R}}
\newcommand{\R}{{\bf R}}
\newcommand{\complex}{{\bf C}}
\newcommand{\intplus}{{\bf N}}
\newcommand{\allint}{{\bf Z}}
\renewcommand{\div}{\mathop{\rm div}}

\newcommand{\e}{{\rm e}}                %% \e = 2.771828...
\renewcommand{\d}{\,{\rm d}}            %% integration measure
\newcommand{\dist}{{\rm dist}}

\newcommand{\cA}{{\cal A}}

\newcommand{\cL}{{\cal L}}
\newcommand{\cM}{{\cal M}}

\newcommand{\cO}{{\cal O}}

\newcommand{\cS}{{\cal S}}
\newcommand{\cT}{{\cal T}}

\newtheorem{theorem}{Theorem}[section]
\newtheorem{lemma}[theorem]{Lemma}

\newtheorem{proposition}[theorem]{Proposition}

\newtheorem{remark}[theorem]{Remark}

\newcommand{\reff}[1]{(\ref{#1})}

\newcommand{\inttwo}{\int_{\real^2}}
\newcommand{\proof}{{\noindent \bf Proof:\ }}
\def\Re{{\rm Re\,}}

\def\pp{{\rm pp}}

\def\div{\mathop{\rm div}}

\def\epsilon{\varepsilon}
\def\phi{\varphi}
\def\weakto{\rightharpoonup}
\def\virgp{\raise 2pt\hbox{,}}
\def\eqdef{\buildrel\hbox{\small{def}}\over =}
\def\build#1_#2^#3{\mathrel{
  \mathop{\kern 0pt#1}\limits_{#2}^{#3}}}
\def\QED{\mbox{}\hfill$\Box$}

\renewcommand{\phi}{\varphi}
\newcommand{\one}{\mathbf{1}}
\newcommand{\spec}{\mathrm{sp}}
\def\bmo{\hbox{\small BMO}}
\def\vmo{\hbox{\small VMO}}
\def\sumetage#1#2{
\sum_{\scriptstyle {#1}\atop\scriptstyle {#2}}
}

\usepackage{amsfonts}
\begin{document}

\title{Uniqueness for the two-dimensional Navier-Stokes equation
   with a measure as initial vorticity}

\author{
Isabelle Gallagher \\
Centre de Math\'ematiques L. Schwartz\\
Ecole polytechnique \\
91128 Palaiseau, France \\
Isabelle.Gallagher@math.polytechnique.fr
\and
Thierry Gallay \\ 
Institut Fourier \\
Universit\'e de Grenoble I \\
38402 Saint-Martin d'H\`eres, France \\
Thierry.Gallay@ujf-grenoble.fr}

\maketitle
\begin{abstract}
We show that any solution of the two-dimensional Navier-Stokes
equation whose vorticity distribution is uniformly bounded in
$L^1(\real^2)$ for positive times is entirely determined by the
trace of the vorticity at $t = 0$, which is a finite measure. When
combined with previous existence results by Cottet, by Giga, Miyakawa,
and Osada, and by Kato, this uniqueness property implies that the
Cauchy problem for the vorticity equation in $\real^2$ is globally
well-posed in the space of finite measures. In particular, this
provides an example of a situation where the Navier-Stokes equation is
well-posed for arbitrary data in a function space that is large enough
to contain the initial data of some self-similar solutions.
\end{abstract}

%%%%%%%%%%%%%%%%%%%%%%%%%%%%%%%%%%%%%%%%%%%%%%%%%%%%%%%%%%%%%%%%%%%%%%%

\section{Introduction}\label{intro} 
\setcounter{equation}{0}

We consider the two-dimensional incompressible Navier-Stokes 
equation
\begin{equation}\label{NS}
  \frac{\partial u}{\partial t} + (u \cdot \nabla) u \,=\,
  \Delta u -\nabla p~, \quad   \div\, u \,=\, 0~, \quad
  x \in \real^2~, \quad t > 0~,
\end{equation}
where $u(x,t) \in \real^2$ denotes the velocity field of the fluid and
$p(x,t) \in \real$ the pressure field. Since this system is very
famous, we do not comment here on its derivation and rather refer to
the monographs \cite{chorin-marsden}, \cite{ladyzhenskaia},
\cite{temam} for a general introduction. The first mathematical
result on the Cauchy problem is due to Leray \cite{leray} who
proved that, for any initial data $u_0 \in L^2(\real^2)$,
system~\reff{NS} has a unique global solution $u \in
C^0([0,+\infty),L^2(\real^2))$ such that~$u(\cdot,0) = u_0$
and~$\nabla u \in L^2((0,+\infty),L^2(\real^2))$. The space
$L^2(\real^2)$ is naturally associated with the Navier-Stokes equation
for two different reasons. First it is the {\sl energy space}, because
the square of the $L^2$ norm of $u$ is the total (kinetic) energy of
the fluid, which is nonincreasing with time. Next, the space
$L^2(\real^2)$ is {\sl scale invariant}, in the sense that $\|\lambda
u_0(\lambda\cdot)\|_{L^2(\real^2)} = \|u_0\|_{L^2(\real^2)}$ for any
$u_0 \in L^2(\real^2)$ and any $\lambda > 0$. This is important
because the transformation
\begin{equation}\label{scaling}
  u(x,t) \,\mapsto\, \lambda u(\lambda x,\lambda^2 t)~, \quad 
  p(x,t) \,\mapsto\, \lambda^2 p(\lambda x,\lambda^2 t)~, \quad 
  \lambda > 0~,
\end{equation}
is a symmetry of \reff{NS}. This invariance was used by Kato
\cite{kato84} to prove that the Navier-Stokes equation in the
$d$-dimensional space $\real^d$ is locally well-posed for arbitrary
data in $L^d(\real^d)$ and even globally well-posed for sufficiently
small data in that space, see also \cite{weissler},
\cite{giga-miyakawa1}. Kato's result was subsequently extended to
larger scale invariant function spaces, such as the homogeneous Besov
space $\dot B^s_{p,q}(\real^d)$ with $s = -1+\frac{d}{p}$ and $p <
\infty$, see Cannone and Planchon~\cite{cannone-book},
\cite{cannone-planchon} and Meyer~\cite{meyer}. A similar analysis was
carried out for the vorticity equation in Morrey spaces by Giga and
Miyakawa \cite{giga-miyakawa2}. One interest of dealing with larger
function spaces is that they may contain initial data which are
homogeneous of degree~$-1$ and therefore give rise to self-similar
solutions of \reff{NS}. This is the case of the Besov space above if
$q = \infty$, or of the larger space $\bmo^{-1}(\real^d)$ introduced
by Koch and Tataru \cite{kochtataru}.  In such spaces, however, it is
not known if the Cauchy problem is well-posed for large data, even
locally in time.

We now return to the two-dimensional case $d=2$ which is simpler for
several reasons. First, the a priori estimates allow in that case to
prove that all solutions are global. For instance, in~\cite{gp},
F.~Planchon and the first author proved that, for arbitrary data in
$\dot B^s_{p,q}(\real^2)$ with~$s = -1+\frac{2}{p}$ and~$p,q <
\infty$, there exists a unique global solution to the Navier-Stokes
equation~\reff{NS}. This result was recently extended by
Germain~\cite{germain} to the larger space~$\vmo^{-1}(\real^2)$,
which is the closure of~$\cS(\real^2)$ in~$\bmo^{-1}(\real^2)$. To our
knowledge, this is the largest space for the velocity field in which
one can solve the Navier-Stokes equation for arbitrary data.  Note
however that~$\vmo^{-1}(\real^2)$ does not contain any non-trivial
homogeneous function of degree~$-1$.

Another specificity of the two-dimensional case is that the 
vorticity~$\omega \eqdef \partial_1 u_2 - \partial_2 u_1$ is a
scalar quantity which satisfies a remarkably simple equation, 
namely
\begin{equation}\label{Veq}
  \frac{\partial \omega}{\partial t} + 
  u \cdot \nabla \omega  \,=\, \Delta \omega~, \quad x \in \real^2~,
  \quad t > 0~.
\end{equation}
The velocity field $u(x,t)$ can be reconstructed from the vorticity
distribution $\omega(x,t)$ by the Biot-Savart law
$$
  u(x,t) \,=\, \frac{1}{2\pi} \inttwo 
  \frac{(x - y)^{\perp}}{|x-y|^2}\, \omega(y,t)\d y\ ,
$$
where $x^\perp = (x_1,x_2)^\perp \eqdef (-x_2,x_1)$. In terms of 
the vorticity, the invariance \reff{scaling} reads
\begin{equation}\label{scaling2}
  \omega(x,t) \,\mapsto\, \lambda^2 \omega(\lambda x,\lambda^2 t)~.
\end{equation}
A natural scale invariant space for the vorticity is thus
$L^1(\real^2)$.  The Cauchy problem for \reff{Veq} in~$L^1(\real^2)$
was studied for instance in~\cite{ben-artzi}, where results analogous
to Leray's and Kato's theorems for the velocity field are 
obtained. However, it is important to realize that a vorticity in 
$L^1(\R^2)$ does not imply a velocity field in $L^2(\R^2)$. Indeed, 
if $u \in L^2(\R^2)$ and if $\omega = \partial_1 u_2 - \partial_2 u_1 
\in L^1(\R^2)$, then it is easy to verify that necessarily
$\inttwo \omega \d x = 0$. Since the integral of $\omega$ (which is
the circulation of the velocity field at infinity) is conserved under
the evolution of \reff{Veq}, it follows that if the initial vorticity
has nonzero integral then the associated velocity field will never be
of finite energy. This ``discrepancy'' between function spaces for the
vorticity and the velocity is in fact specific to the two-dimensional
case. Indeed, if for instance $\omega$ solves the vorticity equation
in $L^{\frac32}(\real^3)$, then the associated velocity field does solve
the Navier-Stokes equation in $L^3(\real^3)$.

In this paper, we study the Cauchy problem for the vorticity equation
\reff{Veq} in~$\cM(\real^2)$,  the space of all finite real measures on 
$\real^2$. If $\mu \in \cM(\real^2)$, the total variation of $\mu$ is 
defined by
$$
   \|\mu\|_{\cM} \,=\, \sup\left\{\inttwo \phi \d \mu \,\Big|\, 
   \phi \in C_0(\real^2)\,,~\|\phi\|_{L^\infty} \le 1\right\}~,
$$
where $C_0(\real^2)$ is the set of all real-valued continuous 
functions on $\real^2$ vanishing at infinity. We recall that
$\cM(\real^2)$ equipped with the total variation norm is a Banach
space, whose norm is invariant under the scaling transformation
\reff{scaling2}. Another useful topology on $\cM(\real^2)$ is the
weak$*$-topology which can be characterized as follows: a sequence
$\{\mu_n\}$ in $\cM(\real^2)$ converges weakly to $\mu$ if $\inttwo
\phi \d\mu_n \to \inttwo \phi \d\mu$ as $n \to \infty$ for all $\phi
\in C_0(\real^2)$. In that case, we write~$\mu_n \weakto \mu$.

Existence of solutions of \reff{Veq} with initial data in 
$\cM(\real^2)$ was first proved by Cottet \cite{cottet},
and independently by Giga, Miyakawa, and Osada \cite{GMO}.
Uniqueness is a more difficult problem. Using a Gronwall-type
argument, it is shown in~\cite{GMO} that uniqueness holds 
if the atomic part of the initial vorticity is sufficiently
small, see also~\cite{kato}. The fact that the size condition only
involves the atomic part of the measure is a consequence of the key
estimate (see~\cite{GMO}) 
$$
  \limsup_{t \to 0} t^{1-\frac1q} \|\e^{t\Delta} \mu \|_{L^q} 
  \,\le\, C_q \|\mu\|_\pp~, \quad 1 <  q \le +\infty~,
$$
where $\|\mu\|_\pp$ denotes the total variation of the atomic
part of $\mu \in \cM(\real^2)$. On the other hand, the case of a 
large Dirac mass was solved recently by C.E.~Wayne and the second 
author \cite{gallay-wayne3} using a completely different approach, 
which we now briefly describe. We first observe that, given any 
$\alpha \in \real$, equation \reff{Veq} has an exact self-similar 
solution given by
\begin{equation}\label{oseendef}
  \omega(x,t) \,=\, \frac{\alpha}{t}\,G\Bigl(
    \frac{x}{\sqrt{t}}\Bigr)~, \quad
  u(x,t) \,=\, \frac{\alpha}{\sqrt{t}}\,v^G\Bigl(
    \frac{x}{\sqrt{t}}\Bigr)~, \quad x \in \real^2~, \quad t > 0~,
\end{equation}
where
\begin{equation}\label{GvGdef}
   G(\xi) \,=\, \frac{1}{4\pi}\,\e^{-|\xi|^2/4}~, \quad
   v^G(\xi) \,=\, \frac{1}{2\pi}\,\frac{\xi^\perp}{|\xi|^2}
   \,\Bigl(1 - \e^{-|\xi|^2/4}\Bigr)~, \quad \xi \in \real^2~.
\end{equation}
This solution is often called the {\sl Lamb-Oseen vortex} with total
circulation $\alpha$. In fact $\omega(x,t)$ is also a solution of the
linear heat equation $\partial_t \omega = \Delta \omega$, because the
nonlinearity in \reff{Veq} vanishes identically due to radial symmetry
(this is again specific to the two-dimensional case). The strategy of
\cite{gallay-wayne3} consists in rewriting \reff{Veq} into
self-similar variables as in \reff{selfdef} below.  Using a pair of
Lyapunov functions, the authors show that the Oseen vortices $\alpha
G$ ($\alpha \in \real$) are the only equilibria of the rescaled
equation. By compactness arguments, they deduce that all solutions
converge in $L^1(\real^2)$ to Oseen vortices as $t \to +\infty$, and
as a byproduct that \reff{oseendef} is the unique solution of
\reff{Veq} such that $\|\omega(\cdot,t)\|_{L^1} \le K$ for all $t > 0$
and $\omega(\cdot,t) \weakto \alpha \delta_0$ as $t \to 0$, 
where~$\delta_0$ is the Dirac mass at the origin.

The goal of the present paper is to solve the uniqueness problem in 
the general case by combining the result of \cite{GMO}, 
which works when the initial measure has small atomic part, 
with the method of \cite{gallay-wayne3}, which allows to handle
large Dirac masses. Our main result is the following:

\begin{theorem}\label{unique}
Let $\mu \in \cM(\real^2)$, and fix $T> 0$, $K > 0$. Then the 
vorticity equation~\reff{Veq} has at most one solution 
$$
  \omega \in C^0((0,T),L^1(\real^2) \cap L^\infty(\real^2))
$$
such that $\|\omega(\cdot,t)\|_{L^1} \le K$ for all $t \in (0,T)$ and 
$\omega(\cdot,t) \weakto \mu$ as $t \to 0$.
\end{theorem}

Here and in the sequel, we say that $\omega(t) \equiv \omega(\cdot,t)$ 
is a (mild) solution of \reff{Veq} on $(0,T)$ if the associated 
integral equation
\begin{equation}\label{Vint}
  \omega(t) \,=\, \e^{(t-t_0)\Delta}\omega(t_0) - \int_{t_0}^t 
  \nabla \cdot \e^{(t-s)\Delta} \Bigl(u(s)\omega(s)\Bigr)\d s
\end{equation}
is satisfied for all $0 < t_0 < t < T$. 

If we combine Theorem~\ref{unique} with the existence results in
\cite{cottet}, \cite{GMO}, \cite{kato}, we conclude that there is a
unique global solution to \reff{Veq} for any initial measure
in~$\cM(\real^2)$. In fact the method we use to prove uniqueness also
implies that this solution depends continuously on the data, so that
the Cauchy problem for the vorticity equation~\reff{Veq} is globally
well-posed in the space $\cM(\real^2)$. If in addition we use the
results in \cite{gallay-wayne3} on the long-time behavior of the
solutions, we obtain the following final statement:

\begin{theorem}\label{final}
For any $\mu \in \cM(\real^2)$, the vorticity equation \reff{Veq} 
has a unique global solution 
$$
  \omega \in C^0((0,\infty),L^1(\real^2) \cap L^\infty(\real^2))
$$
such that $\|\omega(\cdot,t)\|_{L^1} \le \|\mu\|_{\cM}$ for all $t > 0$ 
and $\omega(\cdot,t) \weakto \mu$ as $t \to 0$. This solution depends
continuously on the initial measure $\mu$ in the norm topology of 
$\cM(\real^2)$, uniformly in time on compact intervals. Moreover,
$$
  \inttwo \omega(x,t)\d x \,=\, \alpha \,\eqdef\, \mu(\real^2)~,
  \quad \hbox{for all } t > 0~,
$$
and
\begin{equation}\label{asym}
  \lim_{t \to \infty} t^{1-\frac{1}{p}} \Bigl\|\omega(x,t) - 
  \frac{\alpha}{t}\,G\Bigl(\frac{x}{\sqrt{t}}\Bigr)
  \Bigr\|_{L^p_x} \,=\, 0~, \quad \hbox{for all }
  p \in [1,\infty]~.
\end{equation}
\end{theorem}

Remark that the space $\cM(\real^2)$ does contain nontrivial
homogeneous distributions (the Dirac masses at the origin), 
hence Theorem~\ref{final} gives an example of a situation where the 
Navier-Stokes equation is well-posed for arbitrary data in a function 
space that is large enough to contain the initial data of some 
self-similar solutions (the Oseen vortices). Remark also that 
Oseen's vortex plays a double role in Theorem~\ref{final}: 
it is the unique solution of \reff{Veq} when the initial
vorticity~$\mu$ is a Dirac mass at the origin, and on the other hand it 
describes the long-time behavior of all solutions, see \reff{asym}. 
In fact, it is possible to show that \reff{asym} is a consequence of 
the uniqueness of the solution when $\mu = \alpha \delta_0$, 
see \cite{carpio} and \cite{gigabook}. 

The rest of this paper is devoted to the proof of Theorem~\ref{unique}
and of the continuity statement in Theorem~\ref{final}. 
Before entering the details, let us give a short idea of the argument. 
Previous works on the subject assumed that the initial vorticity 
$\mu$ either has a small atomic part~\cite{GMO},\cite{kato}, or
consists of a single Dirac mass~\cite{gallay-wayne3}. So it is natural 
to decompose~$\mu$ into a finite sum of mutually singular Dirac
masses, and a remainder whose atomic part is arbitrarily small 
(depending on the number of terms in the previous sum). The idea is 
then to use the methods of~\cite{gallay-wayne3} to deal with the 
large Dirac masses, and the argument of~\cite{GMO},\cite{kato} 
to treat the remainder. The difficulty is of course that 
equation~\reff{Veq} is nonlinear so that the interactions between 
the various terms have to be controlled. 

To implement these ideas, we start in Section~\ref{fundamental} by
recalling some general properties of convection-diffusion equations,
of the heat semi-group in self-similar variables, and of the
Biot-Savart law. The proof of Theorem~\ref{unique} begins in
Section~\ref{atomic}, where we decompose the initial measure as
explained above and show that the solution $\omega(x,t)$ also admits a
natural decomposition into a sum of Oseen vortices and a remainder. In
Section~\ref{integral}, we derive the integral equations satisfied by
the remainder terms, and we state a few crucial estimates that will be
proved in an appendix (Section~\ref{estimates}). These results are
used in Section~\ref{gronwall}, where Theorem~\ref{unique} is proved
by a Gronwall-type argument. The same techniques also establish the
continuity claim in Theorem~\ref{final}. 

\medskip\noindent{\bf Notations.}  
We denote by $K_0, K_1, \dots$ our main constants,  
the values of which are fixed throughout the paper. In 
contrast, we denote by $C_0, C_1, \dots$ local constants 
which can take different values in different paragraphs. 
Other positive constants (which are not used anywhere else in
the text) will be generically denoted by $C$. As a general rule, 
we do not distinguish between scalars and vectors in function
spaces: although $u(x,t)$ is a vector field, we write $u \in 
L^2(\real^2)$ and not $u \in L^2(\real^2)^2$. To simplify the
notation, we denote the map $x \mapsto \omega(x,t)$ by $\omega(\cdot,t)$
or just by $\omega(t)$. 

\medskip\noindent{\bf Acknowledgements.} The authors thank 
D.~Iftimie for fruitful discussions. 

%%%%%%%%%%%%%%%%%%%%%%%%%%%%%%%%%%%%%%%%%%%%%%%%%%%%%%%%%%%%%%%%%%%%%%%

\section{Preliminaries}\label{prelim}
\setcounter{equation}{0}

This section is a collection of known results that will be used in
the proof of Theorem~\ref{unique}. 

\subsection{Fundamental solution of a convection-diffusion equation}
\label{fundamental}

We consider the following linear convection-diffusion equation
\begin{equation}\label{Ueq}
  \partial_t \omega(x,t) + U(x,t) \cdot \nabla \omega(x,t) 
  \,=\, \Delta \omega(x,t)~,
\end{equation}
where $x \in \real^2$, $t \in (0,T)$, and $U : \real^2 \times 
(0,T) \to \real^2$ is a (given) time-dependent divergence-free
vector field. The results collected here are due to Osada~\cite{osada},
and to Carlen and Loss~\cite{carlen-loss}. 

Following~\cite{carlen-loss}, we suppose that $U \in C^0((0,T),
L^\infty(\real^2))$ and that
\begin{equation}\label{Ubdd}
  \|U(\cdot,t)\|_{L^\infty(\R^2)} \,\le\, \frac{K_0}{\sqrt t}~,
  \quad 0 < t < T~,
\end{equation}
for some $K_0 > 0$. According to~\cite{osada}, we also assume that
$\Omega \eqdef \partial_1 U_2 - \partial_2 U_1 \in 
C^0((0,T),L^1(\real^2))$ with
\begin{equation}\label{Ombdd}
  \|\Omega(\cdot,t)\|_{L^1(\real^2)} \,\le\, K_0~, \quad 0 < t < T~.
\end{equation}
Then any solution $\omega(x,t)$ of \reff{Ueq} can be
represented as
\begin{equation}\label{omrepr}
  \omega(x,t) \,=\, \inttwo \Gamma_U(x,t;y,s)\omega(y,s)\d y~,
  \quad x \in \real^2~, \quad 0 < s < t < T~,
\end{equation}
where $\Gamma_U$ is the fundamental solution of the 
convection-diffusion equation \reff{Ueq}. The following
properties of $\Gamma_U$ will be useful:

\medskip\noindent{$\bullet$} For any $\beta \in (0,1)$
there exists $K_1 > 0$ (depending only on $K_0$ and $\beta$) 
such that 
\begin{equation}\label{gamm1}
  0 \,<\, \Gamma_U(x,t;y,s) \,\le\, \frac{K_1}{t-s}
  \,\exp\Bigl(-\beta\frac{|x-y|^2}{4(t-s)}\Bigr)\ ,
\end{equation}
for $x,y \in \real^2$ and $0 < s < t < T$, see 
\cite{carlen-loss}. We also have a similar Gaussian lower bound, 
see \cite{osada}. 

\medskip\noindent{$\bullet$} There exists $\gamma \in (0,1)$ 
(depending only on $K_0$) and, for any $\delta > 0$, there exists
$K_2 > 0$ (depending only on $K_0$ and $\delta$) such that 
\begin{equation}\label{gamm2}
  |\Gamma_U(x,t;y,s) - \Gamma_U(x',t';y',s')| \le K_2
  \Bigl(|x{-}x'|^\gamma + |t{-}t'|^{\gamma/2} + 
  |y{-}y'|^\gamma + |s{-}s'|^{\gamma/2}\Bigr)\ ,
\end{equation}
whenever $t-s \ge \delta$ and $t'-s' \ge \delta$, see~\cite{osada}.

\medskip\noindent{$\bullet$} For $0 < s < t < T$ and $x,y \in
\real^2$,  
\begin{equation}\label{gamm3a}
  \inttwo \Gamma_U(x,t;y,s)\d x \,=\, 1~, \quad 
  \inttwo \Gamma_U(x,t;y,s)\d y \,=\, 1~. 
\end{equation}
For $0 < s < r < t < T$ and $x,y \in \real^2$,  
\begin{equation}\label{gamm3b}
  \Gamma_U(x,t;y,s) \,=\, \inttwo \Gamma_U(x,t;z,r) 
  \Gamma_U(z,r;y,s)\d z~.  
\end{equation}

\begin{remark}\label{uptozero}
If $x,y \in \real^2$ and $t > 0$, it follows from~\reff{gamm2} that 
the function $s \mapsto \Gamma_U(x,t;y,s)$ can be continuously extended 
up to $s = 0$, and that this extension (still denoted by $\Gamma_U$) 
satisfies properties~\reff{gamm1} to~\reff{gamm3b}
with $s = 0$.
\end{remark}

\subsection{The heat semiflow in self-similar variables}
\label{heat}

Let $\omega(x,t)$ be a solution of the linear heat equation 
$\partial_t \omega = \Delta \omega$ in $\real^2$. As is well-known, 
it is natural to rewrite this system in terms of the ``self-similar
variables'' $\xi = \frac{x}{\sqrt{t}}$, $\tau = \log(t)$. If we set
\begin{equation}\label{selfdef}
  \omega(x,t) \,=\, \frac1t \,w\Bigl(\frac{x}{\sqrt{t}}\,,\,\log(t)
  \Bigr)~, \quad x \in \real^2~,\quad t > 0~,
\end{equation}
then the new function $w(\xi,\tau)$ is a solution of the rescaled
equation $\partial_\tau w = \cL w$, where $\cL$ is the Fokker-Planck
operator
\begin{equation}\label{Ldef}
  \cL \,=\, \Delta_\xi + \frac12 \xi\cdot\nabla_\xi + 1~.
\end{equation}
This operator is the generator of a $C_0$ semigroup
$S(\tau) = \exp(\tau\cL)$ given by the explicit formula
\begin{equation}\label{Sexplicit} 
  \Bigl(S(\tau)f\Bigr)(\xi) \,=\, \frac{\e^\tau}{4\pi a(\tau)}
  \int_{\real^2} \,\e^{-\frac{|\xi-\xi'|^2}{4a(\tau)}}
  f(\xi'\e^{\frac{\tau}{2}})\d \xi'~, \quad \xi \in \real^2~,
  \quad \tau > 0~,
\end{equation}
where $a(\tau) = 1 - \e^{-\tau}$. The linear operators $\cL$ and 
$S(\tau)$ are studied in detail in (\cite{gallay-wayne1}, Appendix~A). 
For the reader's convenience, we recall here the main properties that 
will be used in the proof of Theorem~\ref{unique}. 

Following \cite{gallay-wayne1}, we introduce for $q \ge 1$ and 
$m \ge 0$ the weighted Lebesgue space $L^q(m)$ defined by
\begin{equation}\label{Lqmdef}
  L^q(m) \,=\, \Bigl\{w \in L^q(\real^2) \,\Big|\, \|w\|_{L^q(m)} 
  < \infty\Bigr\}~, \quad \hbox{where } \|w\|_{L^q(m)} \,=\, 
  \|(1{+}|\xi|^2)^{\frac{m}{2}}w\|_{L^q}~.
\end{equation}
We shall mainly use the Hilbert space $L^2(m)$, which satisfies 
$L^2(m) \hookrightarrow L^1(\real^2)$ if $m > 1$. In this
case, we define the closed subspace
$$
   L^2_0(m) \,=\, \Bigl\{w \in L^2(m) \,\Big|\, \inttwo 
   w(\xi)\d\xi = 0\Bigr\}~, \quad m > 1~.
$$

\begin{proposition}\label{semigroup}
Fix $m > 1$.\\[1mm]
i) There exists $K_3 > 0$ such that, for all $w \in L^2(m)$,
\begin{equation}\label{semi1}
  \|S(\tau) w\|_{L^2(m)} \,\le\, K_3 \|w\|_{L^2(m)}~, 
  \quad \|\nabla S(\tau) w\|_{L^2(m)} \,\le\, 
  \frac{K_3}{a(\tau)^{\frac{1}{2}}}\|w\|_{L^2(m)}~,
\end{equation}
for all $\tau > 0$, where $a(\tau) = 1 - \e^{-\tau}$.\\[1mm]
ii) If moreover $m > 2$ and $w \in L^2_0(m)$, then
\begin{equation}\label{semi2}
  \|S(\tau)w\|_{L^2(m)} \,\le\, K_3 \e^{-\frac{\tau}{2}}\|w\|_{L^2(m)}~,
  \quad \tau \ge 0~.
\end{equation}
iii) More generally, if $q \in [1,2]$ there exists $K_4 > 0$
such that, for all $w \in L^q(m)$,
\begin{equation}\label{semi3}
  \|S(\tau)w\|_{L^2(m)} \,\le\, \frac{K_4}{a(\tau)^{\frac{1}{q} - 
    \frac{1}{2}}}\|w\|_{L^q(m)}~, \quad
  \|\nabla S(\tau)w\|_{L^2(m)} \,\le\,  
   \frac{K_4}{a(\tau)^{\frac{1}{q}}}\|w\|_{L^q(m)}~,
\end{equation}
for all $\tau > 0$. 
\end{proposition} 

\proof The bounds \reff{semi1}, \reff{semi2} are proved in 
(\cite{gallay-wayne1}, Proposition~A.2). Estimate~\reff{semi3} follows 
from \reff{semi1} if we use in addition (\cite{gallay-wayne1}, 
Proposition~A.5). \QED

\medskip
Since the operator $\cL$ has variable coefficients, it does not
commute with spatial derivatives, nor does the associated 
semigroup $S(\tau)$. However, the following useful identity holds:
\begin{equation}\label{Sid}
  \nabla S(\tau) \,=\, \e^{\frac{\tau}2} S(\tau)\nabla~, \quad 
  \tau \ge 0~.
\end{equation}

\subsection{The Biot-Savart law}
\label{Biot-Savart}

Finally we list some basic properties of the Biot-Savart law
\begin{equation}\label{BS2}
  u(x) \,=\, \frac{1}{2\pi} \inttwo \frac{(x - y)^{\perp}}{|x-y|^2}\, 
  \omega(y)\d y~, \quad x \in \real^2~.
\end{equation}
We recall that $L^2(m)$ is the weighted Lebesgue space defined
in \reff{Lqmdef}. 

\begin{proposition}\label{BSprop}
Assume that $\omega \in L^p(\real^2)$ for some $p \in (1,2)$, 
and let $u$ be the vector field defined by \reff{BS2}. Then\\[1mm]
i) $u \in L^q(\real^2)$ where $\frac{1}{q} = \frac{1}{p} - \frac{1}{2}$, 
and there exists $C > 0$ such that 
\begin{equation}\label{HLS}
  \|u\|_{L^q} \,\le\, C \|\omega\|_{L^p}~.
\end{equation}  
ii) $\nabla u \in L^p(\real^2)$ and there exists $C > 0$ such that
\begin{equation}\label{CZ}
  \|\nabla u\|_{L^p} \,\le\, C \|\omega\|_{L^p}~.
\end{equation}  
In addition, $\div(u) = 0$ and $\partial_1 u_2 - \partial_2 u_1 
= \omega$. \\[1mm]
iii) Let $b(x) = (1{+}|x|^2)^{\frac12}$. If $\omega \in L^2(m)$ 
for some $m \in (0,1)$, or $\omega \in L^2_0(m)$ for some
$m \in (1,2)$, then $b^{m-\frac{2}{q}}u \in L^q(\real^2)$ for any
$q \in (2,\infty)$ and there exists $C > 0$ such that
\begin{equation}\label{weightedHLS}
  \|b^{m-\frac{2}{q}} u\|_{L^q} \,\le\, C \|\omega\|_{L^2(m)}~.
\end{equation}  
\end{proposition}

\proof The bound \reff{HLS} is a direct consequence of the 
classical Hardy-Littlewood-Sobolev inequality, see for instance 
(\cite{stein}, Chapter~V, Theorem~1). Estimate~\reff{CZ} holds
because $\nabla u$ is the convolution of $\omega$ with a singular
integral kernel of Calder\'on-Zygmund type, see (\cite{stein}, 
Chapter~II, Theorem~3). Finally, the weighted inequality 
\reff{weightedHLS} is proved in (\cite{gallay-wayne1}, 
Proposition~B.1). \QED

%%%%%%%%%%%%%%%%%%%%%%%%%%%%%%%%%%%%%%%%%%%%%%%%%%%%%%%%%%%%%%%%%%%%%%%

\section{Decomposition of the solution}\label{atomic}
\setcounter{equation}{0}

After these preliminaries, we begin the proof of 
Theorem~\ref{unique}. We fix $\mu \in \cM(\real^2)$, $T> 0$ and~$K > 0$,
and we assume that $\omega \in C^0((0,T),L^1(\real^2) \cap 
L^\infty(\real^2))$ is a solution of the vorticity equation 
\reff{Veq} satisfying $\|\omega(\cdot,t)\|_{L^1} \le K$ for all 
$t \in (0,T)$ and $\omega(\cdot,t) \weakto \mu$ as $t \to 0$.
From~\cite{brezis} we know that~$\omega(x,t)$ coincides for 
$t > 0$ with a classical solution of \reff{Veq} in $\real^2$ 
as constructed for instance in \cite{ben-artzi}. In 
particular $\omega(x,t)$ is smooth for $t > 0$, and 
since the Cauchy problem for~\reff{Veq} is globally well-posed
in $L^1(\real^2) \cap L^\infty(\real^2)$, we could assume without 
loss of generality that~$T = +\infty$. In the sequel, however, 
we keep $T > 0$ arbitrary. 

Since $\mu \in \cM(\real^2)$ is a finite measure, the set 
$E_\pp = \{x\in\real^2\,|\,\mu(\{x\})\neq 0\}$ of all atoms
of $\mu$ is at most countable, and 
$$
   \|\mu\|_\pp \,\eqdef\, \sum_{x \in E_\pp} |\mu(\{x\})| 
   \,\le\, \|\mu\|_{\cM} \,<\, \infty~.
$$
Therefore, given any $\epsilon > 0$, there exists $N \in \intplus$ 
and $z_1, \dots, z_N \in E_\pp$ with $z_i \neq z_j$ for $i \neq j$
such that $\mu$ can be decomposed as
\begin{equation}\label{atomicdec}
  \mu \,=\, \sum_{i=1}^N \alpha_i \delta_{z_i} + \mu_0~, 
\end{equation}
where $\alpha_i = \mu(\{z_i\}) \neq 0$ and $\|\mu_0\|_\pp \le
\epsilon$. Here $\delta_z$ denotes the Dirac mass located at $z \in
\real^2$. Of course, it may happen that $N = 0$ so that $\mu =
\mu_0$, but if the set $E_\pp$ is infinite we have to take~$N$ large
if $\epsilon$ is small. From now on we fix $\epsilon > 0$ and
assume that \reff{atomicdec} holds with $\|\mu_0\|_\pp \le \epsilon$.
We denote
\begin{equation}\label{Mppddef}
  M_{\pp} \,=\, \sum_{i=1}^N |\alpha_i| \,\le\, \|\mu\|_\pp~, 
  \quad\hbox{and}\quad d \,=\, \min\Bigl\{|z_i-z_j| \,\Big|\,
  i,j \in \{1,\dots,N\}\,,~i \neq j\Bigr\}~.
\end{equation}
At the very end of the proof, in Section~\ref{twosolutions}, 
we shall assume that $\epsilon$ is sufficiently small.

Let $u(x,t)$ be the velocity field obtained from $\omega(x,t)$ via
the Biot-Savart law \reff{BS2}. Since for all~$t \in (0,T)$ we 
have~$\|\omega(\cdot,t)\|_{L^1} \le K$, it follows from (\cite{carlen-loss}, 
Theorem~2) that $t^{\frac12}\|u(\cdot,t)\|_{L^\infty} \le CK$ for 
all $t \in (0,T)$, where $C > 0$ is a universal constant. Thus
$\omega(x,t)$ is a solution of the convection-diffusion equation
\reff{Ueq} with $U(x,t) = u(x,t)$, and assumptions \reff{Ubdd}, 
\reff{Ombdd} are satisfied. It follows that $\omega(x,t)$ can
be represented as in \reff{omrepr}, where the fundamental 
solution~$\Gamma_u(x,t;y,s)$ satisfies \reff{gamm1} to \reff{gamm3b}. 
In particular, using Remark~\ref{uptozero}, we have for all~$x \in
\real^2$ and all~$t \in (0,T)$, 
\begin{eqnarray*}
  \omega(x,t) &=& \inttwo \Gamma_u(x,t;y,0)\omega(y,s)\d y \\
  &+& \inttwo (\Gamma_u(x,t;y,s) - \Gamma_u(x,t;y,0))\omega(y,s)\d y~,
  \quad 0 < s < t~.
\end{eqnarray*}
In view of \reff{gamm2}, the second integral in the right-hand side
converges to zero as $s$ goes to zero. On the other hand, since $y
\mapsto \Gamma_u(x,t;y,0)$ is continuous and vanishes at infinity, and
since~$\omega(\cdot,s) \weakto \mu$ as $s \to 0$, we can take the
limit~$s \to 0$ in the first integral and we obtain the following
useful representation:
\begin{equation}\label{representation}
  \omega(x,t) \,=\, \inttwo \Gamma_u(x,t;y,0)\d \mu(y)~, 
  \quad x \in \real^2~, \quad 0 < t < T~.
\end{equation}
Since $\Gamma_u(x,t;y,0)$ is positive and satisfies \reff{gamm3a}, 
it follows that $\|\omega(\cdot,t)\|_{L^1} \le \|\mu\|_{\cM}$ for 
all~$t \in (0,T)$. Thus we can assume that $K = \|\mu\|_{\cM}$ without
loss of generality. 

Inserting \reff{atomicdec} into \reff{representation}, we obtain
the decomposition
\begin{equation}\label{decomp1}
   \omega(x,t) \,=\, \sum_{i=1}^N \omega_i(x,t) + \tilde 
   \omega_0(x,t)~, 
\end{equation}
where
\begin{equation}\label{omegaidef}
  \omega_i(x,t) \,=\, \alpha_i \Gamma_u(x,t;z_i,0)~, \quad 
  x \in \real^2\,,~t \in (0,T)\,,  
\end{equation}
and
\begin{equation}\label{omegazerodef}
  \tilde \omega_0(x,t) \,=\, \inttwo \Gamma_u(x,t;y,0)
  \d \mu_0(y)~, \quad x \in \real^2\,,~t \in (0,T)\,.  
\end{equation}
Thus, although \reff{Veq} is a nonlinear equation, we see that 
the decomposition \reff{atomicdec} of the initial measure induces 
a natural decomposition of the solution $\omega(x,t)$. 
Using the properties of the fundamental solution $\Gamma_u$ 
listed in Section~\ref{fundamental}, one easily obtains
the following results:

\medskip\noindent{$\bullet$}
For all $i \in \{1,\dots,N\}$, $\omega_i \in C^0((0,T),L^1(\real^2) 
\cap L^\infty(\real^2))$ is a solution of \reff{Ueq} with~$U(x,t) 
= u(x,t)$, namely 
\begin{equation}\label{omegaieq}
  \frac{\partial \omega_i}{\partial t} + 
  u \cdot \nabla \omega_i \,=\, \Delta \omega_i~, \quad 
  t \in (0,T)~.
\end{equation}
For any $t \in (0,T)$, $\inttwo \omega_i(x,t)\d x = \alpha_i$, 
$\|\omega_i(\cdot,t)\|_{L^1} = |\alpha_i|$, and
\begin{equation}\label{omegaibdd}
  |\omega_i(x,t)| \,\le\, \frac{K_1 |\alpha_i|}{t} 
  \,\e^{-\beta\frac{|x-z_i|^2}{4t}}~, \quad x \in \real^2~.
\end{equation}
In particular, $\omega_i(\cdot,t) \weakto \alpha_i\delta_{z_i}$ as 
$t \to 0$. 

\medskip\noindent{$\bullet$}
Similarly, $\tilde \omega_0 \in C^0((0,T),L^1(\real^2) \cap 
L^\infty(\real^2))$ is a solution of 
\begin{equation}\label{omegazeroeq}
  \frac{\partial \tilde \omega_0}{\partial t} + 
  u \cdot \nabla \tilde \omega_0 \,=\, \Delta \tilde \omega_0~,
  \quad t \in (0,T)~.
\end{equation}
Moreover, $\|\tilde \omega_0(\cdot,t)\|_{L^1} \le \|\mu_0\|_{\cM} \le 
\|\mu\|_{\cM}$ for all $t \in (0,T)$, and $\tilde \omega_0(\cdot,t) 
\weakto \mu_0$ as $t \to 0$. 

\medskip\noindent
Since $u(x,t)$ is smooth for $t > 0$, it is clear from
\reff{omegaieq}, \reff{omegazeroeq} that $\omega_i(x,t)$ 
and $\tilde \omega_0(x,t)$ are smooth functions of $x \in \real^2$
and $t \in (0,T)$. 

\medskip
For $i \in \{1,\dots,N\}$, we have seen that $\omega_i(x,t)$ 
is a solution of \reff{omegaieq} with a Dirac mass $\alpha_i 
\delta_{z_i}$ as initial data. If we believe in uniqueness, 
we expect that $\omega_i(x,t)$ will be very close, for small 
times, to an Oseen vortex located at $z_i$ with circulation 
$\alpha_i$. Thus if we further decompose
\begin{equation}\label{decomp2}
  \omega_i(x,t) \,=\, \frac{\alpha_i}{t}\,
  G\Bigl(\frac{x-z_i}{\sqrt{t}}\Bigr) +\alpha_i \tilde \omega_i(x,t)~, 
  \quad x \in \real^2~, \quad t \in (0,T)~,  
\end{equation}
where $G$ is defined in \reff{GvGdef}, we expect that the remainder
$\tilde \omega_i(x,t)$ will be small as $t \to 0$. Summarizing, we have
\begin{equation}\label{decomp3}
  \omega(x,t) \,=\, \sum_{i=1}^N \frac{\alpha_i}{t}\, 
    G\Bigl(\frac{x-z_i}{\sqrt{t}}\Bigr) + \tilde \omega(x,t)~, \quad
  u(x,t) \,=\, \sum_{i=1}^N \frac{\alpha_i}{\sqrt{t}}\, 
    v^G\Bigl(\frac{x-z_i}{\sqrt{t}}\Bigr) + \tilde u(x,t)~,
\end{equation}
where
$$
   \tilde \omega(x,t) \,=\,  \tilde \omega_0(x,t) + \sum_{i=1}^N 
\alpha_i \tilde \omega_i(x,t)~, \quad
   \tilde u(x,t) \,=\, \tilde u_0(x,t) + \sum_{i=1}^N  
\alpha_i \tilde u_i(x,t)~,
$$
and where (for $i \in \{0,\dots,N\}$) $\tilde u_i(x,t)$ denotes the 
velocity field associated to $\tilde \omega_i(x,t)$ via the Biot-Savart 
law. In \reff{decomp3}, remark that the explicit terms in the sums 
depend only on the initial measure $\mu$, not on the solution 
$\omega(x,t)$. 

%%%%%%%%%%%%%%%%%%%%%%%%%%%%%%%%%%%%%%%%%%%%%%%%%%%%%%%%%%%%%%%%%%%%%%%

\section{Integral equations and main estimates}\label{integral}
\setcounter{equation}{0}

In this section we derive integral equations for the remainder terms
$\tilde \omega_i(x,t)$ defined in \reff{omegazerodef}
and~\reff{decomp2}, and we also list a few important estimates which
will be proved in Section~\ref{estimates}.  We start in
Section~\ref{diffuse} with~$\tilde \omega_0(x,t)$, which we call the
``diffuse part'' because it is associated to the measure $\mu_0$ which
(by construction) has small or no atomic part. The remaining terms
$\tilde \omega_i(x,t)$ ($i \in \{1,\dots,N\}$), which originate from
the large atoms of the initial measure $\mu$, will be dealt with in
Section~\ref{remainder}.

\subsection{The diffuse part}\label{diffuse}

Let $\tilde \omega_0(x,t)$ be defined by \reff{omegazerodef}. 
Our first result shows that $\tilde \omega_0(x,t)$ is small 
in an appropriate sense as $t \to 0$, because the measure 
$\mu_0$ has a small atomic part. 

\begin{lemma}\label{vanishomegazero}
For any $p \in (1,\infty]$, there exists $K_5 > 0$ (depending 
only on $p$ and $K$) such that
\begin{equation}\label{omegazerolim}
  \limsup_{t \to 0} t^{1-\frac1p} \|\tilde \omega_0(\cdot,t)\|_{L^p} 
  \,\le\, K_5 \|\mu_0\|_\pp~.
\end{equation}  
\end{lemma}

\proof This property is established in (\cite{GMO}, Lemma~4.4) 
in the particular case where~$\tilde \omega_0(\cdot,t) = 
\e^{t\Delta}\mu_0$. By \reff{gamm1}, the fundamental solution
$\Gamma_u(x,t;y,s)$ satisfies a Gaussian upper bound which 
has the same form as the heat kernel $\e^{(t-s)\Delta}(x,y)$, 
so using the same arguments as in \cite{GMO} we immediately 
obtain \reff{omegazerolim}. \QED

\medskip Our next result reflects the fact that $\mu_0(\{z_i\}) = 0$
for $i \in \{1,\dots,N\}$. 

\begin{lemma}\label{estimatetildeomega0}
Assume that $\chi : [0,+\infty) \to \real_+$ is continuous and 
nonincreasing, with $\chi(0) = 1$ and $\chi(r) \to 0$ as 
$r \to \infty$. Then for all $i \in \{1,\dots,N\}$, the following 
estimates hold:
\begin{eqnarray}\label{etildeom}
  \lim_{t \to 0} t^{1-\frac{1}{p}}\|\tilde \omega_0(x,t)
    \chi(\textstyle{\frac{|x-z_i|^2}{t}})\|_{L^p_x} &=& 0~, 
    \quad 1 \le p \le +\infty~, \\ \label{etildeu}
  \lim_{t \to 0} t^{\frac{1}{2}-\frac{1}{q}}\|\tilde u_0(x,t)
    \chi(\textstyle{\frac{|x-z_i|^2}{t}})\|_{L^q_x} &=& 0~, 
    \quad  2 < q \le +\infty~.
\end{eqnarray}
\end{lemma}

\proof See Section~\ref{proof61}. \QED

\medskip
We now derive an integral equation for $\tilde \omega_0(x,t)$. 
Replacing in \reff{omegazeroeq} the velocity field $u(x,t)$ 
with its expression \reff{decomp3} and using Duhamel's formula, 
we obtain for $0 < s < t < T$ the integral representation
\begin{equation}\label{intprelim}
  \tilde \omega_0(t) \,=\, S_N(t,s)\tilde \omega_0(s) - 
  \int_s^t S_N(t,t')(\tilde u(t')\cdot\nabla\tilde \omega_0(t'))
  \d t'~,
\end{equation}
where $\tilde \omega_0(t) \equiv \tilde \omega_0(\cdot,t)$, 
$\tilde u(t) \equiv \tilde u(\cdot,t)$, and $S_N(t,s)$ is the 
evolution operator associated to the convection-diffusion 
equation~\reff{Ueq} with $U(x,t) = \sum_{i=1}^N \frac{\alpha_i}{\sqrt{t}}
\,v^G\Bigl(\frac{x-z_i}{\sqrt{t}}\Bigr)$. 
From Section~\ref{fundamental} we know that
$$
  (S_N(t,s)f)(x) \,=\, \inttwo \Gamma_U(x,t;y,s)f(y)\d y~,
  \quad x \in \real^2~, \quad 0 < s < t~,
$$
where the fundamental solution $\Gamma_U$ satisfies \reff{gamm1} 
to \reff{gamm3b} for some constants $K_1, K_2$ depending on 
$M_\pp$ (but otherwise independent of $N$). By Remark~\ref{uptozero},
$\Gamma_U(x,t;y,s)$ can be continuously extended to~$s = 0$, so that 
$S_N(t,s)$ is well-defined for~$0 \le s < t$. The following properties 
of this operator will be useful: 

\begin{proposition}\label{estimateSN}
Let $p \in [1,\infty]$.\\[1mm]
i) There exists $K_6 > 0$ (depending on $M_{\pp})$
such that, for any measure $\nu \in \cM(\real^2)$,
\begin{equation}\label{estimSN1}
  \|S_N(t,s) \nu\|_{L^p} \,\le\, \frac{K_6}
  {(t-s)^{1-\frac{1}{p}}}\,\|\nu\|_{\cM}~, \quad 0 \le s < t~.
\end{equation}
ii) For any $\gamma \in (0,\frac12)$, there exists $K_7 > 0$ 
(depending on $M_{\pp}$ and $\gamma$) and $t_0 > 0$ (depending
also on $d$) such that, for any function 
$f \in L^1(\real^2)$, 
\begin{equation}\label{estimSN2}
  \|S_N(t,s) \nabla f\|_{L^p} \,\le\, \frac{K_7}{(t-s)^{\frac{3}{2} -
  \frac{1}{p}}}\Bigl(\frac{t}{s}\Bigr)^\gamma \|f\|_{L^1}~, \quad
  0 < s < t < s + t_0~.
\end{equation}
\end{proposition}

\proof See Section~\ref{proof62}. \QED

\begin{remark}
We believe that \reff{estimSN2} holds for $\gamma = 0$ and 
$t_0 = +\infty$, but we were not able to prove that. 
In what follows, we assume without loss of generality that
$t_0 \le T$. 
\end{remark}

As a consequence, if we write $\tilde u(t')\cdot\nabla\tilde
\omega_0(t') = \nabla\cdot(\tilde u(t') \tilde \omega_0(t'))$ in the
right-hand side of \reff{intprelim} and if we use the bound
\reff{estimSN2}, we see that the integral in \reff{intprelim} has a
limit in $L^1(\real^2)$ as $s \to 0$. Moreover, proceeding as in the
proof of \reff{representation}, we obtain $S_N(t,s)\tilde \omega_0(s)
\to S_N(t,0)\mu_0$ as $s \to 0$. Thus $\tilde \omega_0(t)$ satisfies
the integral equation
\begin{equation}\label{omegazeroint}
  \tilde \omega_0(t) \,=\, S_N(t,0)\mu_0 - 
  \int_0^t S_N(t,s)\nabla \cdot(\tilde u(s) \tilde \omega_0(s))
  \d s~, \quad 0 < t < T~.
\end{equation}

\subsection{The atomic part}\label{remainder}

We now fix $i \in \{1,\dots,N\}$ and consider the quantity 
$\tilde \omega_i(x,t)$ defined in \reff{omegaidef}, \reff{decomp2}.
Following~\cite{gallay-wayne1}, \cite{gallay-wayne3}, we introduce the 
self-similar variables
$$
   \xi \,=\, \frac{x-z_i}{\sqrt{t}}~, \quad 
   \tau \,=\, \log(t)~.
$$
We define new functions $\tilde w_i(\xi,\tau)$, $\tilde v_i(\xi,\tau)$ 
by the relations
\begin{equation}\label{tildewvdef}
  \tilde \omega_i(x,t) \,=\, \frac1t \,\tilde w_i\Bigl(
   \frac{x-z_i}{\sqrt{t}}\,,\,\log(t)\Bigr)~, \quad
  \tilde u_i(x,t) \,=\, \frac1{\sqrt{t}} \,\tilde v_i\Bigl(
   \frac{x-z_i}{\sqrt{t}}\,,\,\log(t)\Bigr)~,
\end{equation}
where $x \in \real^2$, $t \in (0,T)$, hence $\xi \in \real^2$, 
$\tau \in (-\infty,\log(T))$. For notational convenience, we 
also define
\begin{equation}\label{tilderel}
   w_i(\xi,\tau) \,=\, \alpha_i G(\xi) + \alpha_i \tilde
   w_i(\xi,\tau)~, \quad v_i(\xi,\tau) \,=\, \alpha_i v^G(\xi) 
   + \alpha_i \tilde v_i(\xi,\tau)~,
\end{equation}
where $G$ and $v^G$ are defined in \reff{GvGdef}. 
In view of \reff{decomp2}, we thus have
\begin{equation}\label{wvdef}
  \omega_i(x,t) \,=\, \frac1t \,w_i\Bigl(
   \frac{x-z_i}{\sqrt{t}}\,,\,\log(t)\Bigr)~, \quad
  u_i(x,t) \,=\, \frac1{\sqrt{t}} \,v_i\Bigl(
   \frac{x-z_i}{\sqrt{t}}\,,\,\log(t)\Bigr)~,
\end{equation}
where $u_i$ is the velocity field associated to $\omega_i$ via 
the Biot-Savart law. 

Inserting these definitions into \reff{omegaieq}, we obtain the 
following evolution equation for $w_i$:
\begin{equation}\label{wieq}
  \frac{\partial w_i}{\partial \tau}(\xi,\tau) + 
  v_i(\xi,\tau)\cdot\nabla w_i(\xi,\tau) + 
  R_i(\xi,\tau)\cdot\nabla w_i(\xi,\tau) \,=\, 
  (\cL w_i)(\xi,\tau)~,
\end{equation}
where $\cL$ is the Fokker-Planck operator \reff{Ldef} and
\begin{equation}\label{Rdef}
  R_i(\xi,\tau) \,=\, \sumetage{j = 1}{j\neq i}^N v_j(\xi - (z_j{-}z_i)
  \e^{-\frac{\tau}{2}},\tau) + \e^{\frac{\tau}{2}}\tilde u_0(\xi 
  \e^{\frac{\tau}{2}}+z_i,\e^\tau)~.
\end{equation}
The corresponding integral equation reads:
\begin{equation}\label{wiinteq1}
  w_i(\tau) \,=\,  S(\tau - \tau_0) w_i(\tau_0) - \int_{\tau_0}^\tau
  S(\tau - \tau') \Bigl(v_i(\tau') + R_i(\tau')\Bigr) 
  \cdot \nabla w_i(\tau') \d \tau'~,
\end{equation}
for $-\infty < \tau_0 < \tau < \log(T)$. Here $S(\tau) = \exp (\tau\cL)$
is the semigroup generated by $\cL$, and $w_i(\tau) \equiv w_i(\cdot,\tau)$, 
$v_i(\tau) \equiv v_i(\cdot,\tau)$. Alternatively, using \reff{Sid}
and the fact that $v_i$, $R_i$ are divergence-free vector fields, 
we have
\begin{equation}\label{wiinteq2}
  w_i(\tau) \,=\,  S(\tau - \tau_0) w_i(\tau_0) - \int_{\tau_0}^\tau
  \e^{-\frac12(\tau-\tau')}\nabla\cdot S(\tau - \tau') 
  \Bigl((v_i(\tau') + R_i(\tau')) w_i(\tau')\Bigr) \d \tau'~.
\end{equation}

It is clear from the definitions that $w_i \in C^0((-\infty,\log(T)),
L^1(\real^2) \cap L^\infty(\real^2))$. Moreover, by~\reff{omegaibdd}, 
we have the pointwise bound
\begin{equation}\label{wibdd}
  |w_i(\xi,\tau)| \,\le\, K_1 |\alpha_i| \,\e^{-\beta |\xi|^2/4}~,
  \quad \xi \in \real^2~, \quad -\infty < \tau < \log(T)~.
\end{equation}
In particular, for any~$m > 1$, the trajectory $\{w_i(\tau)\}$ is 
bounded in the weighted space $L^2(m)$ defined in \reff{Lqmdef}.
Since $w_i(\xi,\tau)$ is smooth, it follows that $w_i \in 
C^0((-\infty,\log(T)),L^2(m))$ for any~$m > 1$. Our next
result shows that $w_i(\tau)$ actually converges to $\alpha_i G$ 
as $\tau \to -\infty$:

\begin{proposition}\label{converge} For any $i \in \{1,\dots,N\}$ 
and any $m > 1$, $w_i(\tau) \to \alpha_i G$ in $L^2(m)$ 
as $\tau \to -\infty$. 
\end{proposition}

\proof See Section~\ref{proof63}. \QED

\medskip This result implies that $\tilde w_i(\tau)$ converges 
to zero in $L^2(m)$ for any $m > 1$. In particular, returning to 
the original variables, we obtain
\begin{equation}\label{vantilde}
  \lim_{t \to 0} t^{1-\frac1p} \|\tilde \omega_i(\cdot,t)\|_{L^p}
  \,=\, 0~, \quad p \in [1,2]~.
\end{equation}

\medskip
We now derive an integral equation for the remainder $\tilde 
w_i(\xi,\tau)$. If we neglect for the moment the term 
$R_i \cdot \nabla w_i$ in \reff{wieq}, and if we replace 
in this equation the functions $w_i, v_i$ by their expressions
\reff{tilderel} and keep only the linear terms in $\tilde w_i, 
\tilde v_i$, we obtain the following equation:
$$
  \frac{\partial \tilde w_i}{\partial \tau} + \alpha_i (v^G  \cdot
  \nabla \tilde w_i +  \tilde v_i \cdot
  \nabla G) \,=\, \cL w_i~.
$$
As is shown in \cite{gallay-wayne3}, this system defines 
a $C_0$ semigroup in $L^2(m)$, which we denote by 
$T_{\alpha_i}(\tau)$. We have the following result, 
which generalizes Proposition~\ref{semigroup}: 

\begin{proposition}\label{talphai}
Fix $\alpha \in \real$ and $m > 1$.\\[1mm]
i) There exists $K_8 > 0$ such that, for all $w \in L^2(m)$,
\begin{equation}\label{talpha1}  
  \|T_{\alpha}(\tau) w\|_{L^2(m)} \,\le\, K_8 \|w\|_{L^2(m)}~, 
  \quad \tau \ge 0~.
\end{equation}
ii) If moreover $m > 2$ and $w \in L^2_0(m)$ then
\begin{equation}\label{talpha2}
  \|T_{\alpha}(\tau) w\|_{L^2(m)} \,\le\, K_8 \,\e^{-\frac{\tau}{2}}
  \|w\|_{L^2(m)}~, \quad \tau \ge 0~.
\end{equation}
iii) Finally if $q \in [1,2]$ and $m > 2$, then $T_{\alpha}(\tau)
\nabla$ can be extended to a bounded operator from~$L^q(m)$
to~$L^2_0(m)$ and there exists $K_9 > 0$ such that
\begin{equation}\label{talpha3}
  \|T_{\alpha}(\tau) \nabla w\|_{L^2(m)} \,\le\, K_9
  \frac{\e^{-\frac{\tau}{2}}}{a(\tau)^{\frac{1}{q}}} 
  \|w\|_{L^q(m)}~, \quad \tau > 0~,
\end{equation}
where $a(\tau) = 1 - \e^{-\tau}$. 
\end{proposition}

\proof See Section~\ref{proof64}. \QED

\begin{remark}
\label{nondecreasingalpha}
One can show that the constants $K_8$, $K_9$ in 
Proposition~\ref{talphai} are uniformly bounded for~$ \alpha$
in compact intervals. 
\end{remark}

Now if we replace in \reff{wieq} the functions $w_i, v_i$ with
their expressions \reff{tilderel} and if we use the above 
notation, we see that $\tilde w_i(\tau)$ is a solution of the
integral equation
$$
  \tilde w_i (\tau) \,=\, T_{\alpha_i}(\tau - \tau_0)\tilde w_i (\tau_0)
  - \int_{\tau_0}^\tau T_{\alpha_i} (\tau - \tau') \Bigl(
  \alpha_i \tilde v_i \cdot \nabla \tilde w_i  + R_i \cdot\nabla (G +
  \tilde w_i)\Bigr)(\tau')\d\tau'~,  
$$
for $-\infty < \tau_0 < \tau < \log(T)$. By Proposition~\ref{converge}, 
$\tilde w_i(\tau) \to 0$ in $L^2(m)$ as $\tau \to -\infty$. 
Thus taking the limit $\tau_0 \rightarrow -\infty$ and using
Proposition~\ref{talphai}, we obtain the desired equation: 
\begin{equation}\label{tildewiint}
  \tilde w_i(\tau) \,=\, -\int_{-\infty}^\tau T_{\alpha_i}(\tau-\tau')
  \nabla \cdot \Bigl(\alpha_i\tilde v_i(\tau')\tilde w_i(\tau') 
  + R_i(\tau') (G + \tilde w_i(\tau'))\Bigr)\d \tau'~.
\end{equation}

%%%%%%%%%%%%%%%%%%%%%%%%%%%%%%%%%%%%%%%%%%%%%%%%%%%%%%%%%%%%%%%%%%%%%%%

\section{The contraction argument}\label{gronwall} 
\setcounter{equation}{0}

This section is devoted to the proof of Theorem~\ref{unique} and
of the continuity statement in Theorem~\ref{final}. By 
\reff{decomp3}, we know that our solution $\omega(x,t)$ 
can be decomposed into a finite sum of Oseen vortices and a 
remainder $\tilde \omega(x,t)$ which is small due to 
\reff{omegazerolim}, \reff{vantilde}. Thus a natural idea 
is to consider the equation satisfied by $\tilde \omega(x,t)$ 
and to apply a Gronwall argument as in \cite{GMO}. However,
this approach requires very precise estimates on the evolution
operator associated to the linearized equation
$$
  {\partial \tilde \omega \over \partial t} + 
  \sum_{i=1}^N \left({\alpha_i \over \sqrt{t}} 
  v^G\Bigl({x-z_i \over \sqrt{t}}\Bigr) \cdot \nabla \tilde
  \omega + \tilde u \cdot \nabla \Bigl({\alpha_i \over t} 
  G\Bigl({x-z_i \over \sqrt{t}}\Bigr)\Bigr)\right) \,=\, 
 \Delta \tilde \omega~,
$$
which are not easy to obtain. Instead we chose to apply a 
Gronwall argument directly to the set of equations \reff{omegazeroint}
\reff{tildewiint}, because the evolution operators 
$S_N(t,s)$ and $T_{\alpha_i}(\tau)$ that appear in these 
equations are simpler to estimate and were already studied 
in \cite{osada}, \cite{gallay-wayne3}. The price to pay with 
this approach is that \reff{tildewiint} still contains some 
linear terms in the right-hand side, which will make the Gronwall 
argument more delicate. 

To make the computations easier to follow, we first deal with 
a single solution in Section~\ref{onesolution}, and in 
Section~\ref{twosolutions} we deduce estimates on the
difference of two solutions which will imply Theorem~\ref{unique}.  
In Section~\ref{continuity} this argument is adapted to prove the 
continuity statement in Theorem~\ref{final}. 

\subsection{Estimates on a single  solution}
\label{onesolution}

Let $\omega$ be a solution of \reff{Veq} satisfying the assumptions
of Theorem~\ref{unique}, and let $u$ be the corresponding velocity 
field. We recall that the initial measure $\mu$ can be decomposed
as in \reff{atomicdec}, with $\|\mu_0\|_\pp \le \epsilon$ for some
$\epsilon > 0$ that will be fixed in Section~\ref{twosolutions}.
According to \reff{decomp3}, $\omega$ and $u$ can be decomposed 
as follows:
$$
 \omega(x,t) \,=\, \sum_{i=1}^N \frac{\alpha_i}{t}\, 
    G\Bigl(\frac{x-z_i}{\sqrt{t}}\Bigr) + \tilde \omega(x,t)~, \quad
  u(x,t) \,=\, \sum_{i=1}^N \frac{\alpha_i}{\sqrt{t}}\, 
    v^G\Bigl(\frac{x-z_i}{\sqrt{t}}\Bigr) + \tilde u(x,t)~,
$$
where
$$
   \tilde \omega(x,t) \,=\,  \tilde \omega_0(x,t) + \sum_{i=1}^N
   \alpha_i \tilde \omega_i(x,t)~, \quad
   \tilde u(x,t) \,=\, \tilde u_0(x,t) + \sum_{i=1}^N \alpha_i   
   \tilde u_i(x,t)~.
$$
Moreover, according to \reff{omegazeroint} and \reff{tildewiint}, 
the remainder terms $\tilde \omega_i$ satisfy the following integral 
equations:

\smallskip\noindent$\bullet$ 
For $i = 0$ and $0 < t < T$, 
\begin{equation}\label{tildew0int2}
  \tilde \omega_0(t) \,=\, S_N(t,0)\mu_0 - \int_0^t S_N(t,s)
  \nabla\cdot(\tilde u(s)\tilde \omega_0(s))\d s~.
\end{equation}

\smallskip\noindent$\bullet$
For $i \in \{1,\dots,N\}$ and $-\infty < \tau < \log(T)$,  
\begin{equation}\label{tildewiint2}
 \tilde w_i(\tau) \,=\, -\int_{-\infty}^\tau T_{\alpha_i}(\tau-\tau')
  \nabla\cdot\Bigl(\alpha_i\tilde v_i(\tau')\tilde w_i(\tau') 
  + R_i(\tau') (G + \tilde w_i(\tau'))\Bigr)\d \tau'~,
\end{equation}
where according to \reff{tildewvdef}
$$
  \tilde w_i(\xi,\tau) \,=\, \e^\tau \tilde 
    \omega_i(\xi\e^{\frac{\tau}{2}},\e^\tau)~, \quad
  \tilde v_i(\xi,\tau) \,=\, \e^{\frac{\tau}{2}} \tilde 
    u_i(\xi\e^{\frac{\tau}{2}},\e^\tau)~.
$$

\medskip\noindent
Fix $m > 2$. For $t \in (0,T)$, we define $M(t) = \max\{M_0(t), 
M_1(t),\dots,M_N(t)\}$, where
$$
   M_0(t) \,=\, \sup_{0 < s \le t} s^{\frac{1}{4}} \|\tilde
   \omega_0(s)\|_{L^{\frac{4}{3}}}~, \quad
   M_i(t) ~=\, \sup_{-\infty < \tau' \le \log(t)} \|\tilde 
   w_i(\tau')\|_{L^2(m)}~, \quad i \in \{1,\dots,N\}~.
$$
We have the following results: 

\begin{proposition}\label{estimateM0}
There exist positive constants $K_{10}, K_{11}$ (depending only on 
$M_{\pp}$) such that
$$
   M_0(t) \,\le\, \delta_1(t) + K_{11} M_0(t) M(t)~, 
   \quad 0 < t < t_0~,
$$
where $t_0 > 0$ is as in Proposition~\ref{estimateSN} and 
$\delta_1(t) \le K_{10}\epsilon$ for $t > 0$ small enough
(depending on~$\mu_0$).
\end{proposition}

\noindent We recall that $M_\pp$ and $d$ are the quantities 
defined in \reff{Mppddef}. 

\medskip\proof
The first term in the right-hand side of \reff{tildew0int2} can be 
estimated as in Lemma~\ref{vanishomegazero}, namely
\begin{equation}\label{muo}
  \limsup_{t \to 0} t^{\frac{1}{4}}\|S_N(t,0)\mu_0\|_{L^{\frac{4}{3}}} 
  \,\le\, K_{10}\|\mu_0\|_{\pp} \,\le\, K_{10} \epsilon~,
\end{equation}
where $K_{10}$ depends only on $M_\pp$. To bound the integral in
\reff{tildew0int2}, we observe that $t^{\frac14} 
\|\tilde \omega_i(t)\|_{L^{\frac43}} \le CM_i(t)$ for $0 < t < T$. 
This is obvious for $i = 0$, whereas for $i \in \{1,\dots,N\}$
we have
$$
  t^{\frac14} \|\tilde \omega_i(t)\|_{L^{\frac43}} \,=\,
  \|\tilde w_i(\log(t))\|_{L^{\frac43}} \,\le\, 
  C \|\tilde w_i(\log(t))\|_{L^2(m)} \,\le\, C M_i(t)~,
$$
since $L^2(m) \hookrightarrow L^{\frac43}(\real^2)$. It follows
that
$$
  t^{\frac14} \|\tilde \omega(t)\|_{L^{\frac43}} \,\le\, 
  M_0(t) + C \sum_{i=1}^N |\alpha_i| M_i(t) \,\le\, C_1 M(t)~,
$$
where $C_1 > 0$ depends only on $M_\pp$. As a consequence, using
\reff{HLS} and H\"older's inequality, we find 
$$
  \|\tilde u(t) \tilde \omega_0(t)\|_{L^1} \,\le\, 
  \|\tilde u(t)\|_{L^4} \|\tilde \omega_0(t)\|_{L^{\frac43}} 
  \,\le\, C \|\tilde \omega(t)\|_{L^{\frac43}} 
  \|\tilde \omega_0(t)\|_{L^{\frac43}} 
  \,\le\, C \frac{M_0(t)M(t)}{t^{\frac12}}~, \quad 0 < t < T~.
$$
Now, using Proposition~\ref{estimateSN}, we obtain for 
$t \in (0,t_0)$
\begin{eqnarray*}
 &&t^{\frac14} \Bigl\|\int_0^t S_N(t,s) \nabla\cdot (\tilde u(s) 
  \tilde\omega_0(s)) \d s \Bigr\|_{L^{\frac43}} 
  \,\le\, t^{\frac{1}{4}} \int_0^t \frac{K_7}{(t-s)^{\frac34}}
  \Bigl(\frac{t}{s}\Bigr)^\gamma \|\tilde u(s) \tilde \omega_0(s)
  \|_{L^1} \d s \\
 &&\qquad \qquad \qquad\qquad \qquad
  \,\le\, t^{\frac{1}{4}} \int_0^t \frac{C}{(t-s)^{\frac34}}
  \Bigl(\frac{t}{s}\Bigr)^\gamma \frac{M_0(s)M(s)}{s^{\frac12}}
  \d s \,\le\, K_{11} M_0(t)M(t)~.
\end{eqnarray*}
Combining this estimate with \reff{muo} we obtain the desired 
result. \QED

\begin{proposition}\label{estimateMi}
There exists a constant $K_{12} > 0$ depending only on $M_\pp$ 
such that, for all~$i \in \{1,\dots,N\}$ and all $t \in (0,T)$, 
$$
  M_i (t) \,\le\, \delta_2(t) + \eta(t) M(t) + K_{12} M_i(t) M(t)~, 
$$
where $\eta(t)$ and $\delta_2(t)$ converge to zero as $t \to 0$. 
\end{proposition}

\proof
We fix $i \in \{1,\dots,N\}$ and estimate successively all terms 
in the right-hand side of \reff{tildewiint2}. Using \reff{Rdef} 
and \reff{tilderel}, we obtain
\begin{equation}\label{Fikdef}
  \int_{-\infty}^\tau T_{\alpha_i}(\tau-\tau')\nabla\cdot
  \Bigl(\alpha_i \tilde v_i(\tau')\tilde w_i(\tau') + R_i(\tau') 
  (G + \tilde w_i(\tau'))\Bigr)\d \tau' \,=\, 
  \sum_{k=1}^6 F_{i,k}(\tau)~,
\end{equation}
where
\begin{eqnarray*}
  F_{i,1}(\tau) &=& \sum_{j\neq i} \int_{-\infty}^\tau T_{\alpha_i} 
    (\tau-\tau') \nabla\cdot \Bigl(\alpha_j v^G (\xi - (z_j{-}z_i) 
    \e^{-\frac{\tau'}{2}}) G\Bigr)\d\tau'~,\\
  F_{i,2}(\tau) &=& \sum_{j\neq i} \int_{-\infty}^\tau 
    T_{\alpha_i} (\tau-\tau') \nabla\cdot \Bigl(\alpha_j \tilde v_j 
    (\xi - (z_j{-}z_i)\e^{-\frac{\tau'}{2}},\tau')G\Bigr)\d\tau'~,\\
  F_{i,3}(\tau) &=& \int_{-\infty}^\tau T_{\alpha_i} (\tau
    -\tau') \nabla\cdot \Bigl(\e^{\frac{\tau'}{2}} \tilde u_0 (\xi
    \e^{\frac{\tau'}{2}} +  z_i,\e^{\tau'}) G\Bigr)\d\tau'~,\\
  F_{i,4}(\tau) &=&  \sum_{j\neq i} \int_{-\infty}^\tau 
    T_{\alpha_i} (\tau-\tau') \nabla\cdot \Bigl(\alpha_j v^G 
    (\xi - (z_j{-}z_i)\e^{-\frac{\tau'}{2}})\tilde w_i(\tau')\Bigr)
    \d\tau'~,\\
  F_{i,5}(\tau) &=& \sum_{j=1}^N \int_{-\infty}^\tau T_{\alpha_i} 
    (\tau-\tau') \nabla\cdot \Bigl(\alpha_j \tilde v_j (\xi - (z_j{-}z_i)
    \e^{-\frac{ \tau'}{2}},\tau')\tilde w_i (\tau')\Bigr)
    \d\tau'~,\\
  F_{i,6}(\tau) &=& \int_{-\infty}^\tau T_{\alpha_i} (\tau
    -\tau') \nabla\cdot \Bigl(\e^{\frac{\tau'}{2}} \tilde u_0 
    (\xi \e^{\frac{\tau'}{2}} +  z_i,\e^{\tau'}) 
    \tilde w_i(\tau')\Bigr)\d\tau'~. 
\end{eqnarray*}
We start with~$F_{i,1} $. Recalling that $\|w\|_{L^2(m)} = 
\|b^m w\|_{L^2}$ with $b(\xi) = (1{+}|\xi|^2)^{\frac12}$, we find
using Proposition~\ref{talphai}
\begin{eqnarray*}
  \|F_{i,1}(\tau)\|_{L^2(m)} &\le&  K_9 \sum_{j\neq i} |\alpha_j| 
    \int_{-\infty}^\tau \frac{\e^{-\frac{\tau - \tau'}{2}}}
    {a(\tau - \tau')^{\frac12}} \|v^G (\xi - (z_j{-}z_i) 
    \e^{-\frac{\tau'}{2}}) b^m G\|_{L^2} \d \tau' \\
  &\le& K_9 \sum_{j\neq i} |\alpha_j| 
    \int_{-\infty}^\tau\frac{\e^{-\frac{\tau - \tau'}{2}}}
    {a(\tau - \tau')^{\frac12}} \,\theta_j(\tau')\|G\|_{L^2(m+1)}
    \d\tau'~,
\end{eqnarray*}
where
$$
  \theta_j(\tau) \,\eqdef\, \sup_{\xi \in \real^2} |b(\xi)^{-1}
  v^G(\xi - (z_j{-}z_i) \e^{-\frac{\tau}{2}})|~, \quad 
  j \neq i~.
$$
Using the explicit expression \reff{GvGdef}, it is easy to verify 
that $\theta_j(\tau) \le C \e^{\frac{\tau}{2}}$ for some $C > 0$ 
depending only on $d$. It follows that
\begin{equation}\label{Fi1}
 \|F_{i,1}(\tau)\|_{L^2(m)} \,\le\, C_1 \e^{\frac{\tau}{2}}~,
 \quad -\infty < \tau < \log(T)~,
\end{equation}
where $C_1 > 0$ depends on $M_\pp$ and $d$. 

To estimate $F_{i,2}$, we write similarly
$$
  \|F_{i,2}(\tau)\|_{L^2(m)} \,\le\, K_9 \sum_{j\neq i} |\alpha_j| 
  \int_{-\infty}^\tau \frac{\e^{-\frac{\tau - \tau'}{2}}}
  {a(\tau - \tau')^{\frac12}} \Bigl\|\tilde v_j(\xi - (z_j{-}z_i)
  \e^{-\frac{\tau'}{2}},\tau') b^m G \Bigr\|_{L^2} \d\tau'~.
$$
Next, we fix $q \in (2,\infty)$ and $\nu \in (0,1)$ such that
$\nu > \frac{2}{q}$. Then, by H\"older's inequality, 
$$
  \Bigl\|\tilde v_j(\xi - (z_j{-}z_i) \e^{-\frac{\tau'}{2}},\tau') 
  b^m G \Bigr\|_{L^2} \,\le\, \Bigl\|b^{\nu - \frac2q} 
  \tilde v_j(\tau')\Bigr\|_{L^q}  \Bigl\|b(\xi - (z_j{-}z_i)
  \e^{-\frac{\tau'}{2}})^{\frac2q-\nu} b^m G\Bigr\|_{L^{\frac{2q}{q-2}}}~.
$$
In view of \reff{weightedHLS} the first factor in the right-hand side 
can be estimated by $C \|\tilde w_j(\tau')\|_{L^2(m)}$, and a direct 
calculation shows that the second one is bounded by 
$C \e^{\tau'(\frac{\nu}{2}-\frac1q)}$, where $C > 0$ depends on $d$. 
Thus
\begin{eqnarray}\nonumber
  \|F_{i,2}(\tau)\|_{L^2(m)} &\le& C \sum_{j\neq i} |\alpha_j| 
    \int_{-\infty}^\tau \frac{\e^{-\frac{\tau - \tau'}{2}}}
    {a(\tau - \tau')^{\frac12}} \|\tilde w_j(\tau')\|_{L^2(m)} 
    \e^{\tau'(\frac{\nu}{2}-\frac{1}{q})}\d\tau' \\ \label{Fi2}
  &\le& C_2 M(\e^\tau) \e^{\tau (\frac{\nu}{2}-\frac1q)}~,
\end{eqnarray}
where $C_2 > 0$ depends on $M_\pp$ and $d$. 

Now we consider $F_{i,3}$. Using Proposition~\ref{talphai} 
and H\"older's inequality, we obtain
\begin{eqnarray}\nonumber
  \|F_{i,3}(\tau)\|_{L^2(m)} &\le&  K_9 \int_{-\infty}^\tau 
  \frac{\e^{-\frac{\tau - \tau'}{2}}}{a(\tau - \tau')^{\frac{1}{2}}} 
  \left\|\e^{\frac{\tau'}{2}} \tilde u_0 (\xi \e^{\frac{\tau'}{2}} + z_i,
  \e^{\tau'}) G^{\frac{1}{2}}\right\|_{L^4} \|b^m G^{\frac{1}{2}}\|_{L^4}
  \d \tau' \\ \label{Fi3}
  &\le& C \int_{-\infty}^\tau \frac{\e^{-\frac{\tau-\tau'}{2}}}
  {a(\tau-\tau')^{\frac12}}\lambda_i(\e^{\tau'})\d \tau' \,\le\, 
  C_3 \lambda_i(\e^{\tau})~,
\end{eqnarray}
where $C_3 > 0$ depends on $M_\pp$ and
$$
  \lambda_i(t) \,=\, \sup_{0<s\le t} s^{\frac{1}{4}} 
  \|\tilde u_0(x,s)\e^{-\frac{|x-z_i|^2}{8s}}\|_{L^4}~.
$$
Applying Lemma~\ref{estimatetildeomega0} with $q = 4$ and 
$\chi(r) = \exp(-r/8)$, we see that $\lambda_i(t)$ converges to zero
as~$t $ goes to~$0$.

For the term $F_{i,4}$, we first remark that
\begin{equation}\label{Atemp}
  \|F_{i,4}(\tau)\|_{L^2(m)} \,\le\,  K_9 \sum_{j\neq i}  |\alpha_j| 
  \int_{-\infty}^\tau \frac{\e^{-\frac{\tau - \tau'}{2}}}
  {a(\tau - \tau')^{\frac12}} \|v^G (\xi - (z_j{-}z_i)  
  \e^{- \frac{\tau'}{2}} ) \tilde w_i(\tau')\|_{L^2(m)}\d \tau'~.
\end{equation}
Since $v^G \in L^\infty(\real^2)$, it follows that
\begin{equation}\label{A}
  \|F_{i,4}(\tau)\|_{L^2(m)} \,\le\, C \sum_{j\neq i}
  |\alpha_j| \int_{-\infty}^\tau \frac{\e^{-\frac{\tau -\tau'}{2}}}
  {a(\tau - \tau')^{\frac12}} \|\tilde w_i(\tau')\|_{L^2(m)} \d \tau'
  \,\le\, C_4' M_i(\e^\tau)~,
\end{equation}
where $C_4' > 0$ depends on $M_\pp$. This bound will be used later on 
when estimating the difference of two solutions. It is not sufficient 
for our present purposes because, unlike in \reff{Fi2}, the prefactor 
of~$M_i(\e^\tau)$ does not converge to zero as $\tau \to -\infty$.
To estimate $F_{i,4}$ more precisely, we observe that, on the one hand,
$$
  |\xi| \le \frac{d}{2} \e^{-\frac{\tau'}{2}} \quad\Rightarrow\quad 
  |v^G (\xi -(z_j{-}z_i) \e^{- \frac{\tau'}{2}})| \le 
  C \e^{\frac{\tau'}{2}}~,
$$
where $C > 0$ depends on $d$. On the other hand the bound \reff{wibdd} 
on $w_i(\xi,\tau)$ implies that
$$
  |\xi| \ge \frac{d}{2} \e^{-\frac{\tau'}{2}} \quad\Rightarrow\quad 
  |\tilde w_i(\xi,\tau')| \le C \e^{-\frac{\beta}{4}|\xi|^2}
  \,\le\, C \e^{-\frac{\beta}{8}|\xi|^2} \e^{-\frac{\beta d^2}{32} 
  \e^{-\tau'}}~.
$$
It follows that
$$
  \|v^G (\xi - (z_j{-}z_i)  \e^{- \frac{\tau'}{2}} )  \tilde w_i
  (\tau')\|_{L^2(m)} \,\le\, C \Bigl(\e^{\frac{\tau'}{2}} 
  \|\tilde w_i(\tau')\|_{L^2(m)} + \zeta(\e^{\tau'})\Bigr)~,
$$
where $C > 0$ depends on $d$ and $\zeta(t) = \exp(-\rho/t)$ for some
$\rho > 0$. Replacing into \reff{Atemp}, we thus find
\begin{equation}\label{Fi4}
  \|F_{i,4}(\tau)\|_{L^2(m)} \,\le\, C_4\Bigl(\e^{\frac{\tau}{2}}
  M_i(\e^\tau) + \zeta (\e^\tau)\Bigr)~,
\end{equation}
where $C_4 > 0$ depends on $M_\pp$ and $d$. 

To estimate~$ F_{i,5}$ we have by Proposition~\ref{talphai}, for 
$1 < p < 2$,
$$
  \|F_{i,5}(\tau)\|_{L^2(m)} \,\le\, K_9 \sum_{j=1}^N 
  |\alpha_j| \int_{-\infty}^\tau \frac{\e^{-\frac{\tau-\tau'}{2}}}
  {a(\tau-\tau')^{\frac{1}{p}}} 
  \|\tilde v_j(\xi - (z_j{-}z_i)\e^{-\frac{\tau'}{2}},\tau') 
  \tilde w_i(\tau')\|_{L^p(m)}\d\tau'~.
$$
If $b(\xi) = (1{+}|\xi|^2)^{\frac12}$, we have using \reff{HLS}
and H\"older's inequality
\begin{eqnarray*}
  \|\tilde v_j \tilde w_i\|_{L^p(m)} \,=\, 
  \|b^m \tilde v_j \tilde w_i\|_{L^p }  & \,\le\, &
  \|b^m \tilde w_i\|_{L^2} \|\tilde v_j\|_{L^{\frac{2p}{2-p}}}  \\
  &\,\le\,& C \|\tilde w_i\|_{L^2(m)} \|\tilde w_j\|_{L^p} 
  \,\le\, C \|\tilde w_i\|_{L^2(m)}\|\tilde w_j\|_{L^2(m)}~.
\end{eqnarray*}
It follows that
\begin{equation}\label{Fi5}
  \|F_{i,5} (\tau) \|_{L^2(m)} \,\le\, C \sum_{j=1}^N  
  |\alpha_j| M_j(\e^\tau) M_i(\e^\tau) \,\le\, 
  C_5 M_i(\e^\tau) M(\e^\tau)~,
\end{equation}
where $C_5 > 0$ depends on $M_{\pp}$.

Finally we consider the last term, $F_{i,6}$. Choosing $p = 
\frac43 \in (1,2)$, we obtain as above
\begin{eqnarray*}
  \|F_{i,6}(\tau)\|_{L^2(m)} &\le& K_9 
    \int_{-\infty}^\tau \frac{\e^{-\frac{\tau - \tau'}{2}}}
    {a(\tau - \tau')^{\frac{1}{p}}} \Bigl\|\e^{\frac{\tau'}{2}} 
    \tilde u_0 (\xi \e^{\frac{\tau'}{2}}  + z_i, \e^{\tau'})
    \tilde w_i(\tau')\Bigr\|_{L^p(m)}\d \tau' \\ 
  &\le &  K_9 \int_{-\infty}^\tau \frac{\e^{-\frac{\tau - \tau'}{2}}}
    {a(\tau - \tau')^{\frac{1}{p}}} \Bigl\|\e^{\frac{\tau'}{2}} 
    \tilde u_0 (\xi \e^{\frac{\tau'}{2}}  + z_i,
    \e^{\tau'}) \Bigr\|_{L^{\frac{2p}{2-p}}} \|\tilde w_i
    (\tau')\|_{L^{2}(m)}\d \tau'~.
\end{eqnarray*}
Using \reff{HLS}, we find (with $t = \e^\tau$)
$$
  \|\e^{\frac{\tau}{2}} \tilde u_0 (\xi \e^{\frac{\tau}{2}}  + z_i,
  \e^\tau)\|_{L^{\frac{2p}{2-p}}} \,=\, t^{1-\frac{1}{p}} 
  \|\tilde u_0 (\cdot,t)\|_{L^{\frac{2p}{2-p}}}
  \,\le\, C t^{1-\frac{1}{p}} \|\tilde \omega_0 (\cdot,t)\|_{L^p} 
  \,\le\, C M_0(t)~,
$$
hence finally
\begin{equation}\label{Fi6}
  \|F_{i,6}(\tau)\|_{L^2(m)} \,\le\, C_6 M_0(\e^\tau) M_i
  (\e^\tau)~,
\end{equation}
where $C_6 > 0$ depends on $M_\pp$. Collecting estimates \reff{Fi1} 
to \reff{Fi6}, we obtain the desired bound on $M_i(t)$. This concludes 
the proof of Proposition~\ref{estimateMi}. 
\QED

\medskip
Note that Propositions~\ref{estimateM0} and \ref{estimateMi}
together imply that
\begin{equation}\label{estimateallMi}
  M(t) \,\le\, \delta(t) + \eta(t) M(t) + K_{13}M(t)^2~, 
  \quad 0 < t < t_0~,
\end{equation}
where $K_{13} > 0$ depends only on $M_\pp$, $\eta(t)$ goes to zero 
as $t \to 0$, and $\delta(t) \le K_{10} \epsilon$ if $t > 0$ is 
small enough. Both functions $\eta(t)$, $\delta(t)$ depend on the
full initial measure $\mu$, not only on $M_\pp$ and $d$. 

\subsection{The uniqueness proof}\label{twosolutions}

This section is devoted to the end of the proof of Theorem~\ref{unique}.
Let $\omega^{(1)}$ and $\omega^{(2)}$ be two solutions of \reff{Veq}
satisfying the assumptions of Theorem~\ref{unique} with the same 
initial measure $\mu$. Each solution can be decomposed as in 
\reff{decomp3}, namely
$$
  \omega^{(\ell)}(x,t) \,=\, \sum_{i=1}^N \frac{\alpha_i}{t}\, 
  G\Bigl(\frac{x-z_i}{\sqrt{t}}\Bigr) + \tilde \omega^{(\ell)}(x,t)~, 
  \quad \tilde \omega^{(\ell)}(x,t) \,=\, \tilde \omega_0^{(\ell)}(x,t) + 
  \sum_{i=1}^N \alpha_i \tilde \omega_i^{(\ell)} (x,t)~, 
$$
for $\ell \in \{1,2\}$. Estimate \reff{estimateallMi} becomes, with 
obvious notations, 
\begin{equation}\label{Mk}
  M^{(\ell)}(t) \,\le\, \delta(t) + \eta(t) M^{(\ell)}(t)
  + K_{13} M^{(\ell)}(t)^2~, \quad \ell \in \{1,2\}~, \quad 
  0 < t < t_0~.
\end{equation}
Now, we define $\Delta(t) = \max\{\Delta_0(t),\Delta_1(t),\dots,
\Delta_N(t)\}$, where
$$
  \Delta_0(t) \,=\, \sup_{0 < s \le t} s^{\frac{1}{4}} \|\tilde
  \omega_0^{(1)}(s) - \tilde\omega_0^{(2)}(s)\|_{L^{\frac43}}~, 
$$
and
\begin{equation}\label{Deltaidef}
  \Delta_i(t) \,=\,  \sup_{-\infty < \tau' \le \log(t)}  
  \|\tilde w_i^{(1)}(\tau') - \tilde w_i^{(2)}(\tau') \|_{L^2(m)}~,
  \quad i \in \{1,\dots,N\}~.
\end{equation}
Here and in the sequel, $\tilde w_i^{(\ell)}(\xi,\tau) \,=\, \e^\tau 
\tilde \omega_i^{(\ell)}(\xi\e^{\frac{\tau}{2}},\e^\tau)$ for $\ell \in 
\{1,2\}$. We have the following result:

\begin{proposition}\label{estimatedelta}
There exists a constant $K_{14} > 0$ depending only on $M_{\pp}$ 
such that
$$
   \Delta(t) \,\le\, \eta(t) \Delta(t) + K_{14}\Bigl(
   M^{(1)}(t) + M^{(2)}(t)\Bigr)\Delta(t) + \zeta(t)~,
   \quad 0 < t < t_0~,
$$
where $\eta(t)$ goes to zero as $t \to 0$ and $\zeta(t) =
C \e^{-\rho/t}$ for some $\rho > 0$. Moreover, 
$$
   \Delta(t) \,\le\, \eta(t) \Delta(t) + K_{14}\Bigl(
   M^{(1)}(t) + M^{(2)}(t)\Bigr)\Delta(t) + K_{14}
   \int_0^t \frac{\Delta(s)}{(t-s)^{\frac12}s^{\frac12}} 
   \d s~.
$$
\end{proposition}

\proof
The argument consists in mimicking the proofs of 
Propositions~\ref{estimateM0} and~\ref{estimateMi} above. 
We start by estimating $\Delta_0(t)$. We have of course
$$
  \tilde \omega_0^{(1)}(t) - \tilde \omega_0^{(2)}(t) \,=\,
   -\int_0^t S_N(t,s)\nabla\cdot\Bigl(\tilde u^{(1)}(s) \tilde 
  \omega^{(1)}_0(s) - \tilde u^{(2)}(s)\tilde 
  \omega^{(2)}_0(s)\Bigr)\d s~.
$$
If we write $\tilde u^{(1)}\tilde \omega^{(1)}_0 - \tilde u^{(2)}
\tilde \omega^{(2)}_0 = (\tilde u^{(1)} - \tilde u^{(2)}) \tilde 
\omega^{(1)}_0 + \tilde u^{(2)}(\tilde \omega^{(1)}_0 - 
\tilde \omega^{(2)}_0)$ and if we proceed exactly as in the proof 
of Proposition~\ref{estimateM0}, we obtain
$$
  \Delta_0(t) \,\le\, C_0 \Delta(t) \Bigl(M^{(1)}(t) + M^{(2)}(t)
  \Bigr)~, \quad 0 < t < t_0~,
$$
where $C_0 > 0$ depends only on $M_{\pp}$.

We now bound $\Delta_i (t)$ for $i \in \{1,\dots,N\}$. Let 
$G_{i,k}(\tau) = F_{i,k}^{(1)}(\tau) - F_{i,k}^{(2)}(\tau)$
for $k \in \{1,\dots,6\}$, where $F_{i,k}^{(1)}$ and $F_{i,k}^{(2)}$
are defined in analogy with \reff{Fikdef}. Then obviously
$G_{i,1} = 0$. The quadratic terms $G_{i,k}$ for $k \in 
\{5,6\}$ can be estimated as in the case of $\Delta_0$ above. 
In view of \reff{Fi5}, \reff{Fi6}, we thus find
$$
  \|G_{i,5}(\tau)\|_{L^2(m)} 
   + \|G_{i,6}(\tau)\|_{L^2(m)} \,\le\, C_1\Delta(\e^\tau) 
   \Bigl(M^{(1)}(\e^\tau)+ M^{(2)}(\e^\tau)\Bigr)~, 
$$
for $\tau \in (-\infty,\log(T))$, where $C_1 > 0$ depends on $M_\pp$.
It remains to bound the linear terms $G_{i,k}$ for $k \in \{2,3,4\}$.
Proceeding as in the proofs of \reff{Fi2}, \reff{Fi4}, we obtain
$$
   \|G_{i,2}(\tau)\|_{L^2(m)} \,\le\, C_2 \Delta(\e^\tau) 
   \e^{\tau (\frac{\nu}{2}-\frac1q)}~, \quad 
   \|G_{i,4}(\tau)\|_{L^2(m)} \,\le\, C_4\Bigl(\e^{\frac{\tau}{2}}
   \Delta_i(\e^\tau) + \zeta(\e^\tau)\Bigr)~,
$$
where $C_2, C_4$ depend on $M_\pp$ and $d$, and $\zeta(t) = C 
\e^{-\rho/t}$ for some $\rho > 0$. Furthermore, using the analogue 
of \reff{A}, we have
$$
   \|G_{i,4}(\tau)\|_{L^2(m)} \,\le\, C_4' \int_{-\infty}^\tau 
   \frac{\e^{-\frac{\tau - \tau'}{2}}}{a(\tau - \tau')^{\frac12}}
   \|\tilde w_i^{(1)}(\tau') - \tilde w_i^{(2)}(\tau')\|_{L^2(m)}  
   \d \tau'~,
$$
where $C_4'$ depends on $M_\pp$. Returning to the original time 
variable $t = \e^\tau$, we thus find
$$
   \|G_{i,4}(\log(t))\|_{L^2(m)} \,\le\, C_4' \int_0^t 
   \frac{\Delta_i(s)}{(t-s)^{\frac12}s^{\frac12}}\d s~, \quad 
   0 < t < T~.
$$
Finally, according to \reff{Fi3}, we have the bound
\begin{eqnarray*}
  \|G_{i,3}(\tau)\|_{L^2(m)} &\le& C \int_{-\infty}^\tau 
  \frac{\e^{-\frac{\tau - \tau'}{2}}}{a(\tau - \tau')^{\frac12}} 
  \Bigl\|\e^{\frac{\tau'}{2}} \tilde u_0^{(1)}(\xi \e^{\frac{\tau'}{2}} 
  + z_i,\e^{\tau'}) - \e^{\frac{\tau'}{2}} \tilde u_0^{(2)} 
  (\xi \e^{\frac{\tau'}{2}} + z_i,\e^{\tau'}) \Bigr\|_{L^4} 
  \d \tau' \\
  &\le& C \int_{-\infty}^\tau \frac{\e^{-\frac{\tau-\tau'}{2}}}
  {a(\tau-\tau')^{\frac12}} \Delta_0(\e^{\tau'})\d \tau' 
  \,\le\, C_3 \Delta_0(\e^{\tau})~,
\end{eqnarray*}
which is sufficient for our purposes since $\Delta_0(t) \le C_0 
\Delta(t)(M^{(1)}(t) + M^{(2)}(t))$. Collecting all these 
estimates, we obtain the desired bounds on $\Delta(t)$. 
This concludes the proof of Proposition~\ref{estimatedelta}. \QED

\medskip\noindent
{\bf Proof of Theorem~\ref{unique}:}  
Let $\tilde K = \max\{K_5,K_{10}\}$, where $K_5$ is as in 
Lemma~\ref{vanishomegazero} and $K_{10}$ as in 
Proposition~\ref{estimateM0}. Assume that $\epsilon > 0$ is 
sufficiently small so that
\begin{equation}\label{howsmall}  
   16 K_{13} \tilde K \epsilon \,\le\, 1~, \quad 
   \hbox{and}\quad 16 K_{14} \tilde K \epsilon \,\le\, 1~,
\end{equation}
where $K_{13}$ is as in \reff{Mk} and $K_{14}$ as in 
Proposition~\ref{estimatedelta}. Finally, choose $t_1 \in (0,t_0]$
sufficiently small so that
$$
  \eta(t) \,\le\, \frac14~, \quad \hbox{and} \quad 
  \delta(t) \,\le\, \tilde K \epsilon~, \quad \hbox{for }
  0 < t \le t_1~,
$$
where $\delta(t), \eta(t)$ are as in \reff{Mk} and $\eta(t)$ 
appears in Proposition~\ref{estimatedelta} as well. 
We shall prove that $\Delta(t) = 0$ for $ \in (0,t_1]$, hence 
$\omega^{(1)}(t) = \omega^{(2)}(t)$ on this time interval. 
Since~$\omega^{(\ell)}(t_1) \in L^1(\real^2) \cap L^\infty(\real^2)$
and since the Cauchy problem is well-posed in that space, 
it will follow that $\omega^{(1)}(t) = \omega^{(2)}(t)$ for all 
$t \in (0,T)$. 

We claim that $M^{(\ell)}(t) \le 2\tilde K \epsilon$ for $\ell \in 
\{1,2\}$ and $t \in (0,t_1]$. Indeed, by Lemma~\ref{vanishomegazero}
and Proposition~\ref{converge}, this is true at least for $t > 0$ 
sufficiently small. On the other hand, it follows from \reff{Mk} 
that
$$ 
   M^{(\ell)}(t) \,\le\, \tilde K \epsilon + \frac14 M^{(\ell)}(t)
   + K_{13} M^{(\ell)}(t)^2~, \quad 0 < t \le t_1~,
$$
hence $M^{(\ell)}(t) < 2\tilde K \epsilon$ as long as $K_{13}
M^{(\ell)}(t) < \frac14$. Since $K_{13}(2\tilde K\epsilon) \le \frac18$, 
this proves the claim. 

Now, it follows from \reff{howsmall} and 
Proposition~\ref{estimatedelta} that 
$$
   \Delta(t) \,\le\, \frac12 \Delta(t) + \zeta(t)~, \quad
    \hbox{and}\quad 
   \Delta(t) \,\le\, \frac12 \Delta(t) +  K_{14}
   \int_0^t \frac{\Delta(s)}{(t-s)^{\frac12}s^{\frac12}} 
   \d s~,
$$
for $t \in (0,t_1]$. The first inequality implies that 
$\Delta(t) = \cO(t^\infty)$ as $t \to 0$. In view of 
Lemma~\ref{decroittrop} below, the second bound then implies that
$\Delta(t) = 0$ for $t \in (0,t_1]$. This concludes the proof 
of Theorem~\ref{unique}. \QED

\begin{lemma}\label{decroittrop}
Let $f : [0,T] \to \real_+$ be a continuous function satisfying
$$
  f(t) \,\le\, K \int_0^t \frac{f(s)}{(t-s)^{\frac12}s^{\frac12}}\d s~,
  \quad 0 \le t \le T~,
$$
for some $K > 0$. If $f(t) = \cO(t^\alpha)$ as $t \to 0$ for all 
$\alpha > 0$, then $f \equiv 0$. 
\end{lemma}

\proof Given $\alpha \ge 0$, we define
$$
  F_\alpha(t) \,=\, \sup_{0 <  s \le t} \frac{f(s)}{s^\alpha}~,
  \quad 0 < t \le T~.
$$
If $F_\alpha(T) < \infty$, we have for $t \in (0,T]$:
$$
  \frac{f(t)}{t^\alpha} \,\le\, \frac{K}{t^\alpha} 
  \int_0^t \frac{s^\alpha F_\alpha(t)}{(t-s)^{\frac12}s^{\frac12}}
  \d s \,=\, K B(\textstyle{\frac12},\alpha{+}\textstyle{\frac12}) 
  F_\alpha(t)~,
$$
where
$$
  B(p,q) \,=\, \int_0^1 (1{-}x)^{p-1} x^{q-1} \d x \,=\, 
  \frac{\Gamma(p)\Gamma(q)}{\Gamma(p+q)}~, \quad p,q > 0~.
$$
It follows that $F_\alpha(T) \le K B(\frac12,\alpha{+}\frac12)
F_\alpha(T)$. Now, if $f(t) = \cO(t^\infty)$ we can take 
$\alpha > 0$ large enough so that $K B(\frac12,\alpha{+}\frac12) < 1$.
Then $F_\alpha(T) = 0$, which implies $f \equiv 0$. \QED

\subsection{The continuity proof}\label{continuity}

In this section we prove that the (unique) solution of \reff{Veq}
depends continuously on the initial data in the norm topology 
of~$\cM(\R^2)$, as stated in Theorem~\ref{final}. The arguments are 
very similar to those leading to the uniqueness theorem, 
except for the fact that the initial measures associated to both 
solutions are now different. So we shall merely sketch the proof 
and emphasize where the arguments of the previous sections must 
be adapted to infer continuity.

Fix $\mu^{(1)} \in \cM(\R^2)$, and assume that $\mu^{(2)}$ is 
another finite measure satisfying $\|\mu^{(1)} - \mu^{(2)}\|_{\cM}
\le \delta$ for some sufficiently small $\delta > 0$. This implies
in particular that the large atoms of $\mu^{(1)}, \mu^{(2)}$ are 
located at the same points in $\real^2$. More precisely, we can assume 
that both measures are decomposed as in Section~\ref{atomic},
namely
$$
   \mu^{(\ell)} \,=\, \sum_{i=1}^N \alpha_i^{(\ell)} \delta_{z_i} 
   + \mu_0^{(\ell)}~, \quad \ell \in \{1,2\}~,
$$ 
where $\alpha_i^{(\ell)} = \mu^{(\ell)}(\{z_i\}) \neq 0$ and  
$\|\mu_0^{(\ell)}\|_\pp \le \epsilon$. The parameter 
$\epsilon > 0$ is independent of $\delta$ and will be assumed 
to satisfy a smallness condition similar to \reff{howsmall}.
By construction, we have
$$
  \|\mu^{(1)} - \mu^{(2)}\|_{\cM} \,=\, \sum_{i=1}^N |\alpha_i^{(1)}
  -\alpha_i^{(2)}| + \|\mu_0^{(1)} - \mu_0^{(2)}\|_{\cM} 
  \,\le\, \delta~.
$$
For $\ell \in \{1,2\}$, let $\omega^{(\ell)} \in C^0((0,+\infty),
L^1(\real^2) \cap L^\infty(\real^2))$ be the solution of 
\reff{Veq} with initial data $\mu^{(\ell)}$. Each solution can be 
decomposed as in \reff{decomp3}, namely
$$
  \omega^{(\ell)}(x,t) \,=\, \sum_{i=1}^N \frac{\alpha_i^{(\ell)}}{t}\, 
  G\Bigl(\frac{x-z_i}{\sqrt{t}}\Bigr) + \tilde \omega^{(\ell)}(x,t)~, 
  \quad \tilde \omega^{(\ell)}(x,t) \,=\, \tilde \omega_0^{(\ell)}(x,t) + 
  \sum_{i=1}^N \alpha_i^{(\ell)} \tilde \omega_i^{(\ell)} (x,t)~.
$$
Using the same notations as in the previous section, our goal is 
to control the quantity $\Delta(t) = \max\{\Delta_0(t),\Delta_1(t),
\dots,\Delta_N(t)\}$, where
$$
  \Delta_0(t) \,=\, \sup_{0 < s \le t} \|\tilde
  \omega_0^{(1)}(s) - \tilde\omega_0^{(2)}(s)\|_{L^1} +  
  \sup_{0 < s \le t} s^{\frac{1}{4}} \|\tilde \omega_0^{(1)}(s) 
  - \tilde\omega_0^{(2)}(s)\|_{L^{\frac43}}~, 
$$
and $\Delta_i(t)$ is defined by \reff{Deltaidef} for $i \in 
\{1,\dots,N\}$. We shall prove that, for any $\nu \in (0,1)$, 
there exists $T > 0$ and $C > 0$ (both independent of $\delta$)
such that $\Delta(T) \le C\delta^\nu$. In particular, this 
implies
\begin{equation}\label{contdef}
   \sup_{0 < t \le T} \|\omega^{(1)}(t) - 
   \omega^{(2)}(t)\|_{L^1} \,\le\, C \delta^\nu~.
\end{equation}
Since the Cauchy problem for \reff{Veq} in $L^1(\real^2)$ is 
globally well-posed and since the solution is a locally Lipschitz 
function of the initial data in that space, uniformly in 
time on compact intervals, it follows that \reff{contdef} holds 
for any $T > 0$. This proves the continuity claim in 
Theorem~\ref{final}. 

\medskip
To bound $\Delta_0(t)$ we write
\begin{eqnarray*}
  \tilde \omega_0^{(1)} (t) - \tilde \omega_0^{(2)} (t) & \,=\, & 
  \Bigl(S_N^{(1)}(t,0) - S_N^{(2)}(t,0)\Bigr)
  \mu^{(1)}_0  + S_N^{(2)}(t,0) \Bigl(\mu^{(1)}_0 - \mu^{(2)}_0 \Bigr) \\
  &-& \int_0^t \Bigl(S_N^{(1)}(t,s) - S_N^{(2)}(t,s)\Bigr) 
    \nabla \cdot(\tilde u^{(1)}(s) \tilde \omega^{(1)}_0(s))\d s \\
  &-& \int_0^t  S_N^{(2)}(t,s)  \nabla \cdot\Bigl(\tilde u^{(1)}(s) 
    \tilde \omega^{(1)}_0(s) - \tilde u^{(2)}(s) \tilde
    \omega^{(2)}_0(s)\Bigr)~, 
\end{eqnarray*}
where~$S_N^{(\ell)}(t,s)$ denotes the evolution operator
associated to the convection--diffusion equation~\reff{Ueq} with 
$U(x,t) = \sum_{i=1}^N \frac{\alpha_i^{(\ell)}}{\sqrt{t}}
\,v^G\Bigl(\frac{x-z_i}{\sqrt{t}}\Bigr)$. The difference 
$S_N^{(1)}(t,s)-S_N^{(2)}(t,s)$ is estimated using the following
variant of Proposition~\ref{estimateSN}, which can be proved by a 
standard perturbation argument (we omit the details). 

\begin{proposition}\label{estimatediffSN}
Let $p \in [1,\infty]$ and let~$S_N^{(1)}$ and~$S_N^{(2)}$ be
defined as above.\\[1mm]
i) There exists $K_{15} > 0$  independent of~$ \delta $
such that, for any measure $\nu \in \cM(\real^2)$,
$$
  \Bigl\|\Bigl(S_N^{(1)}(t,s) - S_N^{(2)}(t,s) \Bigr)\nu\Bigr\|_{L^p} 
  \,\le\, \frac{K_{15}\delta}{(t-s)^{1-\frac{1}{p}}}\,\|\nu\|_{\cM}~, 
  \quad 0 \le s < t~.
$$
ii) For any $\gamma \in (0,\frac12)$, there exists $K_{16} > 0$ and 
$t_0 > 0$ (both independent of~$ \delta $) such that, for any 
function $f \in L^1(\real^2)$, 
$$
  \Bigl\|\Bigl(S_N^{(1)}(t,s) - S_N^{(2)}(t,s) \Bigr) \nabla f\Bigr\|_{L^p}
  \,\le\, \frac{K_{16} \delta}{(t-s)^{\frac{3}{2} -
  \frac{1}{p}}}\Bigl(\frac{t}{s}\Bigr)^\gamma \|f\|_{L^1}~, \quad
  0 < s < t < s + t_0~.
$$
\end{proposition}

Using Propositions~\ref{estimateSN} and~\ref{estimatediffSN}, and 
proceeding as in Section~\ref{twosolutions} above, it is not 
difficult to show that
$$
  \Delta_0(t) \,\le\, C_0 \delta (1 + M^{(1)}(t)^2) 
  + C_0\Bigl(M^{(1)}(t) + M^{(2)}(t)\Bigr)\Delta(t)~, 
  \quad 0 < t < t_0~,
$$
where $C_0 > 0$ is independent of $\delta$ and $t_0 > 0$ is as 
in Proposition~\ref{estimatediffSN}. 

To bound $\Delta_i(t)$ for~$i \in \{1,\dots,N\}$, we consider 
the integral equations of the form \reff{tildewiint2} satisfied 
by the rescaled functions $\tilde w_i^{(1)}(\tau)$ and 
$\tilde w_i^{(2)}(\tau)$, and we estimate the difference of 
both expressions. To this end, we clearly need a bound on the 
linear operator $(T_{\alpha}(\tau) - T_{\alpha'}(\tau))\nabla$
with $\alpha \neq \alpha'$. This is the content of the following 
proposition, whose proof is again left to the reader.

\begin{proposition}\label{difftalphai}
Fix $\alpha \in \real$, $m > 2$, $\nu \in (0,\frac12)$ and~$q
\in [1,2]$. Then there exists $K_{17} > 0$ and~$ \beta_0 > 0$ such that 
if~$|\beta| \le \beta_0$ then, for all $f \in L^q(m)$, 
$$
  \Bigl\|\Bigl(T_{\alpha+\beta}(\tau) -  T_{\alpha}(\tau)
  \Bigr)\nabla f\Bigr\|_{L^2(m)} \,\le\, K_{17}|\beta|
  \frac{\e^{- \nu \tau}}{a(\tau)^{\frac{1}{q}}}
  \|f\|_{L^q(m)}~, \quad \tau > 0~,
$$
where $a(\tau) = 1 - \e^{-\tau}$. 
\end{proposition}

As in Section~\ref{twosolutions}, we define $G_{i,k}(\tau) = 
F_{i,k}^{(1)}(\tau) - F_{i,k}^{(2)}(\tau)$ for $k \in \{1,\dots,6\}$,
where $F_{i,k}^{(1)}$ and $F_{i,k}^{(2)}$ are defined in analogy 
with \reff{Fikdef}. Arguing as in the previous sections, we find 
the following estimates:
\begin{eqnarray*}
  \|G_{i,1}(\tau)\|_{L^2(m)} &\le& C_1 \delta \,\e^{\nu \tau} ~,\\
  \|G_{i,2}(\tau)\|_{L^2(m)} &\le& C_2 \delta \,\e^{\nu \tau} 
    M^{(1)}(\e^\tau)+ C_2 \e^{\nu \tau} \Delta(\e^\tau)~,\\
  \|G_{i,3}(\tau)\|_{L^2(m)} &\le& C_3 \delta M^{(1)}(\e^\tau) + 
    C_3 \Delta_0(\e^\tau)~,\\
  \|G_{i,4}(\tau)\|_{L^2(m)} &\le& C_4 \delta M^{(1)}(\e^\tau)  
    + C_4 \Bigl(\e^{\nu\tau} \Delta(\e^\tau) + \zeta(\e^\tau)\Bigr)~,
\end{eqnarray*}
where $\nu \in (0,\frac12)$ and $\zeta(t) = C \e^{-\rho/t}$ for some 
$\rho > 0$. Here and in the sequel, all constants are independent 
of $\delta$. As in Section~\ref{twosolutions}, the term~$G_{i,4}(\tau)$ 
can also be estimated as follows:
$$
  \|G_{i,4}(\log(t))\|_{L^2(m)} \,\le\, C_4 \delta M^{(1)}(t)  
  + C_4 \int_0^t \frac{\Delta_i(s)}{(t-s)^{\frac12}s^{\frac12}}
  \d s~.
$$
Finally
$$
  \sum_{k \in \{5,6\}} \|G_{i,k}(\tau)\|_{L^2(m)} 
  \,\le\, C_5 \delta M^{(1)}(\e^\tau)^2 + C_5 \Delta(\e^\tau) 
  \Bigl(M^{(1)}(\e^\tau) + M^{(2)}(\e^\tau)\Bigr)~.
$$
Combining these estimates with the above bound on $\Delta_0(t)$, 
we finally obtain
\begin{equation}\label{bddfin1}
  \Delta (t) \,\le\, K_{18} \delta (1 + M^{(1)}(t)^2) 
  + \eta(t) \Delta(t) + K_{19} \Bigl(M^{(1)}(t) +  M^{(2)}(t)\Bigr) 
  \Delta (t) + \zeta(t)~,
\end{equation}
as well as
\begin{eqnarray}\nonumber
  \Delta (t) &\le& K_{18} \delta (1 + M^{(1)}(t)^2)
  + \eta(t) \Delta(t) + K_{19} \Bigl(M^{(1)}(t) +  M^{(2)}(t)\Bigr) 
  \Delta (t) \\ \label{bddfin2}
  &+& K_{20} \int_0^t \frac{\Delta(s)}{ 
  (t-s)^{\frac{1}{2}} s^{\frac{1}{2}}}\d s~.
\end{eqnarray}
Now, proceeding as in the proof of Theorem~\ref{unique} in
Section~\ref{twosolutions}, we can choose $\epsilon > 0$ sufficiently
small and then $t_1 \in (0,t_0]$ sufficiently small (both $\epsilon$
and $t_1$ independent of $\delta$) so that $\eta(t) \le \frac14$ and
$K_{19}(M^{(1)}(t) + M^{(2)}(t)) \le \frac14$ for all $t \in (0,t_1]$.
The bounds \reff{bddfin1}, \reff{bddfin2} then imply that, for any
$\nu \in (0,1)$, there exists $K_{21} > 0$ (independent of $\delta$)
such that $\Delta(t) \le K_{21}\delta^\nu$ for all $t \in (0,t_1]$,
which is the desired result. Indeed, we have the following lemma,
which is a generalization of Lemma~\ref{decroittrop}:

\begin{lemma}\label{decroittropdelta}
Let $f : [0,T] \to \real_+$ be a continuous function satisfying
$$
  f(t) \,\le\, C_1 \delta + C_2 \int_0^t \frac{f(s)}{(t-s)^{\frac12}
  s^{\frac12}}\d s~, \quad 0 \le t \le T~,
$$
for some $C_1,C_2 > 0$ and some $\delta \in [0,1]$. Suppose moreover
that $f(t) \le C_1 \delta +  \zeta (t)$ for all $t \in [0,T]$, 
where $\zeta(t) = C_3 \exp(-\rho/t)$ for some $\rho > 0$. 
Then for any $\nu \in (0,1)$, there exists a constant $C_4 > 0$
(independent of $\delta$) such that $f(t) \le C_4 \delta^\nu$ for 
all $t \in [0,T]$. 
\end{lemma}

\begin{remark}
What the proof really shows is that $\Delta(t_1) \le C\delta 
\log(C/\delta)^\gamma$ for some large $\gamma > 0$. Thus
our method fails to show that the solution of \reff{Veq} is
a locally Lipschitz function of the initial data. This is 
because we chose to apply the Gronwall argument directly
to equation \reff{tildewiint}, see the discussion at the beginning 
of Section~\ref{gronwall}. 
\end{remark}

%%%%%%%%%%%%%%%%%%%%%%%%%%%%%%%%%%%%%%%%%%%%%%%%%%%%%%%%%%%%%%%%%%%%%%%

\section{Appendix}\label{estimates}
\setcounter{equation}{0}

In this final section we prove the main estimates stated in 
Section~\ref{integral}.

\subsection{Proof of Lemma~\ref{estimatetildeomega0}}
\label{proof61} 

We argue as in the proof of (\cite{GMO}, Lemma~4.4).  
Fix $p \in [1,\infty]$, $q \in (2,\infty]$, and $i \in \{1,\dots,N\}$. 
Without loss of generality, we can assume that $z_i = 0$. 
Using the definition \reff{omegazerodef}, the bound~\reff{gamm1}, and
the properties of the Biot-Savart law, it is easy to show that there
exists $C_1 > 0$ such that
\begin{equation}\label{tdbdd}  
  \|\tilde \omega_0(\cdot,t)\|_{L^p} \,\le\, \frac{C_1}{t^{1-\frac1p}}
  \|\mu_0\|~, \quad 
  \|\tilde u_0(\cdot,t)\|_{L^q} \,\le\, \frac{C_1}{t^{\frac12-\frac1q}}
  \|\mu_0\|~, \quad t > 0~.
\end{equation}
Fix any $\delta > 0$. Since $\mu_0(\{0\}) = 0$ by assumption, 
there exists $r > 0$ such that $|\mu_0|(B_{4r})
\le \delta/C_1$, where $B_r = \{x\in\real^2\,|\, |x| \le r\}$
and $|\mu_0|$ denotes the total variation measure associated 
with $\mu_0$. We decompose
$$
 \tilde \omega_0(x,t) \,=\, \int_{B_{4r}} \Gamma_u(x,t;y,0)
 \d \mu_0(y) +  \int_{\real^2 \setminus B_{4r}} 
 \Gamma_u(x,t;y,0) \d \mu_0(y) \,\eqdef\, 
 \omega^{(1)}(x,t) + \omega^{(2)}(x,t)~.
$$
We also have $\tilde u_0(x,t) = u^{(1)}(x,t) + u^{(2)}(x,t)$, 
where $u^{(j)}$ is the velocity field obtained from 
$\omega^{(j)}$ via the Biot-Savart law \reff{BS2}. By construction, 
\begin{equation}\label{ombd1}
  \sup_{t > 0} t^{1-\frac1p} \|\omega^{(1)}(\cdot,t)\|_{L^p} \,\le\,  
    C_1 |\mu_0|(B_{4r}) \,\le\, \delta~, \quad
  \sup_{t > 0} t^{\frac12-\frac1q} \|u^{(1)}(\cdot,t)\|_{L^q} \,\le\, 
    C_1 |\mu_0|(B_{4r}) \,\le\, \delta~.
\end{equation}
To bound $\omega^{(2)}(x,t)$, we further decompose
$$
  \omega^{(2)}(x,t) \,=\, \omega^{(2)}(x,t)\one_{B_{2r}}(x)
  + \omega^{(2)}(x,t)\one_{\real^2\setminus B_{2r}}(x)
  \,\eqdef\, \omega^{(3)}(x,t) + \omega^{(4)}(x,t)~.
$$
Accordingly we set $u^{(2)}(x,t) = u^{(3)}(x,t) + u^{(4)}(x,t)$. 
Using \reff{gamm1}, we find
$$
   |\omega^{(3)}(x,t)| \,\le\, \one_{B_{2r}}(x)
   \int_{\real^2 \setminus B_{4r}} \frac{K_1}{t}
   \,\e^{-\beta\frac{|x-y|^2}{4t}}\d|\mu_0|(y) \,\le\, 
   \,\e^{-\beta\frac{r^2}{2t}} \inttwo \frac{K_1}{t} 
   \,\e^{-\beta\frac{|x-y|^2}{8t}}\d|\mu_0|(y)~,
$$
since $|x-y| \ge 2r$ in the first integral. It follows that
\begin{equation}\label{ombd2}
   t^{1-\frac1p} \|\omega^{(3)}(\cdot,t)\|_{L^p} +  
   t^{\frac12-\frac1q} \|u^{(3)}(\cdot,t)\|_{L^q} \,\le\, C 
   \,\e^{-\beta\frac{r^2}{2t}}\|\mu_0\|_{\cM} 
   \,\build\hbox to 10mm{\rightarrowfill}_{t \to 0}^{}\, 0~.
\end{equation}
Finally $\|\omega^{(4)}(x,t)\chi(|x|^2/t)\|_{L^p_x} \le
\chi(4r^2/t)\|\omega^{(4)}(\cdot,t)\|_{L^p}$, hence
\begin{equation}\label{ombd3}  
   t^{1-\frac1p} \|\omega^{(4)}(x,t)\chi(|x|^2/t)\|_{L^p_x} 
   \,\le\, C_1 \chi(4r^2/t)\|\mu_0\|_{\cM} 
   \,\build\hbox to 10mm{\rightarrowfill}_{t \to 0}^{}\, 0~.
\end{equation}
Combining \reff{ombd1}, \reff{ombd2}, \reff{ombd3}, we obtain
$\limsup_{t\to 0} t^{1-\frac1p} \|\tilde \omega_0(x,t)
\chi(|x|^2/t)\|_{L^p_x} \le \delta$. Since $\delta > 0$ 
was arbitrary, this proves \reff{etildeom}.

To bound $u^{(4)}(x,t)$, we use yet another decomposition:
$$
  u^{(4)}(x,t) \,=\, u^{(4)}(x,t)\one_{B_r}(x)
  + u^{(4)}(x,t)\one_{\real^2\setminus B_r}(x)
  \,\eqdef\, u^{(5)}(x,t) + u^{(6)}(x,t)~.
$$
Using the Biot-Savart law \reff{BS2}, we find
$$
  |u^{(5)}(x,t)| \,\le\, \one_{B_r}(x) 
  \int_{\real^2\setminus B_{2r}} \frac{C}{|x-y|}
  |\omega^{(4)}(y,t)|\d y \,\le\, \frac{C}r
  \one_{B_r}(x) \|\omega^{(4)}(\cdot,t)\|_{L^1}
  \,\le\, \frac{C}r \one_{B_r}(x) 
  \|\mu_0\|~,
$$
hence $t^{\frac12-\frac1q} \|u^{(5)}(\cdot,t)\|_{L^q} \le 
C (t/r^2)^{\frac12-\frac1q} \|\mu_0\|_{\cM} \to 0$ as
$t \to 0$. Finally, 
$$
  t^{\frac12-\frac1q} \|u^{(6)}(x,t)\chi(|x|^2/t)\|_{L^q_x}
  \,\le\, \chi(r^2/t) t^{\frac12-\frac1q}
  \|u^{(4)}(\cdot,t)\|_{L^q} \,\le\, C \chi(r^2/t)
  \|\mu_0\|_{\cM} \,\build\hbox to 10mm{\rightarrowfill}_{t \to
  0}^{}\, 0~.
$$
Summarizing, we have shown $\limsup_{t\to 0} t^{\frac12-\frac1q} 
\|\tilde u_0(x,t)\chi(|x|^2/t)\|_{L^q_x} \le \delta$, which
implies \reff{etildeu}. \QED

\subsection{Proof of Proposition~\ref{estimateSN}}
\label{proof62}

Estimate \reff{estimSN1} follows immediately from the bound 
\reff{gamm1} on the integral kernel $\Gamma_U(x,t;y,s)$. To
prove \reff{estimSN2}, we first remark that it is sufficient
to establish this estimate for $p = 1$. Indeed, once this is done, 
we obtain using \reff{estimSN1}:
\begin{eqnarray*}
&& \|S_N(t,s)\nabla f\|_{L^p} \,=\, \Bigl\|S_N\Bigl(t,\frac{t{+}s}2\Bigr) 
   S_N\Bigl(\frac{t{+}s}2,s\Bigr)\nabla f\Bigr\|_{L^p}
   \,\le\, K_6 \Bigl(\frac2{t{-}s}\Bigr)^{1-\frac1p}
   \Bigl\|S_N\Bigl(\frac{t{+}s}2,s\Bigr)\nabla f\Bigr\|_{L^1} \\
&&\qquad \,\le\, K_6 \Bigl(\frac2{t{-}s}\Bigr)^{1-\frac1p} 
   K_7 \Bigl(\frac2{t{-}s}\Bigr)^{\frac12}
   \Bigl(\frac{1{+}t/s}2\Bigr)^{\gamma}\|f\|_{L^1}
   \,\le\, K_6 K_7 \Bigl(\frac2{t{-}s}\Bigr)^{\frac32-\frac1p}
  \Bigl(\frac{t}{s}\Bigr)^{\gamma}\|f\|_{L^1}~.
\end{eqnarray*}
It remains to prove \reff{estimSN2} for $p = 1$. We proceed in 
two steps: 

\bigskip\noindent{\bf Step 1 : the case of 1 vortex}\\[1mm]
Fix $\alpha \in \real$, and let $S_1(t,s)$ be the evolution 
operator associated to the non-autonomous equation
$$
  \frac{\partial \omega}{\partial t}(x,t) + 
  \frac{\alpha}{\sqrt{t}}\,v^G\Bigl(\frac{x}{\sqrt{t}}\Bigr)
  \cdot \nabla \omega (x,t) \,=\, \Delta \omega(x,t)~,
  \quad x \in \real^2~, \quad t > 0~.
$$
Due to the particular form of the convection term, it is
natural to rewrite this equation in the self-similar 
variables $\xi = x/\sqrt{t}$, $\tau = \log(t)$. Defining
$w(\xi,\tau)$ as in \reff{selfdef}, we obtain the equivalent 
equation
\begin{equation}\label{S1equiv}
  \frac{\partial w}{\partial \tau}(\xi,\tau) + 
  \alpha v^G(\xi)\cdot\nabla w(\xi,\tau) \,=\, 
  (\cL w)(\xi,\tau)~,
\end{equation}
where $\cL$ is given by \reff{Ldef}. We shall show that the
autonomous equation \reff{S1equiv} defines a strongly continuous
semigroup in $L^1(\real^2)$, which we denote by $\cS_1(\tau)$. We
claim that, for any~$\tau > 0$, the operator~$\cS_1(\tau)\nabla$ can
be extended to a bounded operator on~$L^1(\real^2)$. Moreover, for 
any $\gamma > 0$, there exists $C_1 > 0$ (depending on $\gamma$ and 
$|\alpha|$) such that, for any $w \in L^1(\real^2)$,
\begin{equation}\label{S1bdd}
  \|\cS_1(\tau)\nabla w\|_{L^1} \,\le\, C_1\frac{\e^{-(\frac12-\gamma)\tau}}
  {a(\tau)^{\frac12}} \|w\|_{L^1}~, \quad \tau > 0~,
\end{equation}
where $a(\tau) = 1 - \e^{-\tau}$. If we return to the original 
variables, we see that \reff{S1bdd} is equivalent to~\reff{estimSN2}
with $N = 1$, $p = 1$, and $t_0 = +\infty$. 

To prove \reff{S1bdd}, we introduce the Banach space $X \hookrightarrow
L^1(\real^2)$ defined by
$$
   X \,=\, \Bigl\{w \in L^1(\real^2) \,\Big|\, w = \partial_1 f_1 + 
   \partial_2 f_2 \hbox{ with } f_1, f_2 \in L^1(\real^2)\Bigr\}~,
$$
equipped with the norm
$$
  \|w\|_X \,=\, \|w\|_{L^1} + 
  \inf\Bigl\{\|f_1\|_{L^1} + \|f_2\|_{L^1} \,\Big|\, w = \partial_1 f_1 + 
   \partial_2 f_2\Bigr\}~.
$$
We also consider the auxiliary equation for the vector field
$f = (f_1,f_2)$:
\begin{equation}\label{S1aux}
  \frac{\partial f}{\partial \tau} + 
  \alpha v^G \div f \,=\, \Bigl(\cL - \frac12\Bigr)f~.
\end{equation}
Using a fixed point argument as in the proof of Lemma~\ref{auxT} 
below, it is straightforward to show that \reff{S1aux} defines a 
strongly continuous semigroup in $L^1(\real^2)^2$, which we denote 
by $\cT_1(\tau)$. Moreover, there exists $\tau_0 > 0$ and $C_2 > 0$ 
such that, for all $f \in L^1(\real^2)^2$,
\begin{equation}\label{cT1bdd}
  \|\cT_1(\tau)f\|_{L^1} \,\le\, C_2\|f\|_{L^1}~, \quad
  \|\nabla \cT_1(\tau)f\|_{L^1} \,\le\, \frac{C_2}{a(\tau)^{\frac12}}
  \|f\|_{L^1}~, \quad 0 < \tau \le \tau_0~.
\end{equation}
The evolutions defined by \reff{S1equiv} and \reff{S1aux} are 
related via
$$
   \div(\cT_1(\tau)f) \,=\, \cS_1(\tau)\div f~, \quad f \in X~, 
   \quad \tau > 0~.
$$
This shows that, for $\tau \in (0,\tau_0]$, $\cS_1(\tau)\nabla$ can 
be extended to a bounded operator from $L^1(\real^2)$ into $X$ with 
bound $C_2 a(\tau)^{-\frac12}$; in particular \reff{S1bdd} holds for 
$\tau \in (0,\tau_0]$. Moreover $\cS_1(\tau)$ is a strongly continuous 
semigroup in $X$. Thus, to prove \reff{S1bdd} for all times, it 
remains to show that, for any $\gamma > 0$, there exists $C_3 > 0$ 
such that $\|\cS_1(\tau)\|_{\cL(X)} \le C_3 \,\e^{-(\frac12-\gamma)\tau}$. 
Equivalently, we shall show that the spectral radius of $\cS_1(\tau)$ 
in $X$ satisfies $\rho_\spec(\cS_1(\tau)) \le \e^{-\frac{\tau}{2}}$ for 
all $\tau \ge 0$. 

To prove this, we argue exactly as in (\cite{gallay-wayne3},
Sections~4.1 and 4.2). We first observe that~$\cS_1(\tau)$ is a 
compact perturbation of $S(\tau) = \exp(\tau\cL)$, and it is easy 
to verify using \reff{Sexplicit} that the spectral radius of 
$S(\tau)$ in $X$ satisfies $\rho_\spec(S(\tau)) = \e^{-\frac{\tau}{2}}$
for all $\tau \ge 0$. Thus, it remains to show that all eigenvalues
of the generator $L = \cL - \alpha v^G\cdot\nabla$ of $\cS_1(\tau)$ are
contained in the half-plane~$\{\lambda \in \complex \,|\, \Re(\lambda)
\le -\frac12\}$.  Assume on the contrary that some $\lambda \in \complex$
with $\Re(\lambda) > -\frac12$ is an eigenvalue of $L$ in $X$. 
Since $L$ is rotation invariant, we can use polar coordinates in
$\real^2$ and assume that the eigenfunction $\phi$ associated to
$\lambda$ has the form $\phi(r\cos\theta,r\sin\theta) = \psi(r)
\,\e^{in\theta}$ for some $n \in \allint$. If we study the
differential equation satisfied by $\psi$, we find as in 
(\cite{gallay-wayne3}, Lemma~4.5) that
$$
  \psi(r) \,\sim\, A r^{2\lambda-2} + B r^{-2\lambda} 
  \,\e^{-r^2/4}~, \quad r \to +\infty~,
$$
for some $A,B \in \complex$. Now, since $\phi \in X$ and 
$\Re(\lambda) > -\frac12$, we must have $A = 0$, hence $\phi$ 
has Gaussian decay at infinity, and $\inttwo \phi(\xi)\d\xi = 0$. 
In particular, $\phi$ lies in the Hilbert space
$$
   Y_0 \,=\,\Bigl\{w \in L^2(\real^2,\complex) \,\Big|\, 
   \inttwo G^{-1}|w(\xi)|^2\d\xi < \infty\,,~\inttwo w(\xi)\d\xi
   = 0\Bigr\}~.
$$
But it is proved in \cite{gallay-wayne3} that $L$ is self-adjoint 
in $Y_0$ with spectrum $\{-\frac{n}{2}\,|\, n \in \intplus\,,~n \ge 1\}$ 
and~$v^G\cdot\nabla$ is skew-symmetric in the same space $Y_0$. 
Thus we necessarily have $\Re(\lambda) \le -\frac12$, which is a 
contradiction. 

\bigskip\noindent{\bf Step 2 : the case of N vortices}\\[1mm]
We now assume that $N \ge 2$ and we study the evolution 
operator $S_N(t,s)$ associated to the equation
\begin{equation}\label{SNequation} 
  \frac{\partial \omega}{\partial t}(x,t) + \sum_{i=1}^N 
  \frac{\alpha_i}{\sqrt{t}}\,v^G\Bigl(\frac{x-z_i}{\sqrt{t}}\Bigr)
  \cdot \nabla \omega (x,t) \,=\, \Delta \omega(x,t)~.
\end{equation}
As we shall see, if $t/d^2 \ll 1$ where $d = \min\{|z_i-z_j|\,|\, 
i \neq j\}$, the $N$ convection terms in \reff{SNequation} 
are nearly decoupled, and we can bound $S_N(t,s)$ using the previous
estimates on $S_1(t,s)$. 

Let~$\chi : \real^2 \to [0,1]$ be a smooth function equal to one
for~$|x| \le \frac14$ and zero for~$|x| \ge 
\frac13$. For~$i \in \{1,\dots,N\}$ we set~$\chi_i(x) = \chi((x{-}z_i)/d)$
and we define~$\chi_0$ such that~$\sum_{i=0}^N \chi_i(x) = 1$
for all $x \in \real^2$. Observe that $0 \le \chi_0 \le 1$ and that
there exists~$C_1 > 0$ (independent of~$d$) such that~$\|\sum_{i=0}^N|
\nabla \chi_i|\,\|_{L^\infty} \le C_1 d^{-1}$ and~$\|\sum_{i=0}^N|
\Delta \chi_i|\,\|_{L^\infty} \le C_1 d^{-2}$ for all~$i \in \{0,\dots,N\}$.

If $\omega(x,t)$ satisfies \reff{SNequation}, then 
for all $i \in \{0,\dots,N\}$ the function $\omega_i(x,t) \eqdef
\chi_i(x)\omega(x,t)$ is a solution of
$$
  \frac{\partial \omega_i}{\partial t} + \frac{\alpha_i}{\sqrt{t}}
  \,v^G\Bigl(\frac{x-z_i}{\sqrt{t}}\Bigr)\cdot\nabla \omega_i 
  \,=\, \Delta \omega_i -\div (R_i\omega) + Q_i\omega~,
$$
where $\alpha_0 = 0$ and 
\begin{eqnarray*}
  R_i(x,t) &=& \sum_{j\neq i} \frac{\alpha_j}{\sqrt{t}}
    \,v^G\Bigl(\frac{x-z_j}{\sqrt{t}}\Bigr)\chi_i(x) 
    + 2\nabla\chi_i(x)~, \\
  Q_i(x,t) &=& \sum_{j=1}^N \frac{\alpha_j}{\sqrt{t}}
    \,v^G\Bigl(\frac{x-z_j}{\sqrt{t}}\Bigr)\cdot\nabla\chi_i(x) 
    + \Delta\chi_i(x)~.
\end{eqnarray*}
By construction, $R_i(x,t)$ and $Q_i(x,t)$ are smooth functions 
of $x \in \real^2$ and $t \ge 0$. Moreover, if~$R(x,t) = 
\sum_{i=0}^N |R_i(x,t)|$ and $Q(x,t) = \sum_{i=0}^N |Q_i(x,t)|$, 
there exists $C_2 > 0$ (independent of $d$) such that
\begin{equation}\label{RQbdd}
   \|R(\cdot,t)\|_{L^\infty} \,\le\, \frac{C_2}{d}~, \quad
   \|Q(\cdot,t)\|_{L^\infty} \,\le\, \frac{C_2}{d^2}~, \quad
   t \ge 0~.
\end{equation}
If we denote by $\tilde S_i(t,s)$ the evolution operator associated to 
the $i$th vortex, we find the following integral equation
\begin{equation}\label{SNint}
  \omega_i(t) \,=\, \tilde S_i(t,s)\omega_i(s) + \int_s^t 
  \tilde S_i(t,t') \Bigl(-\div (R_i(t')\omega(t')) + 
  Q_i(t')\omega(t')\Bigr)\d t'~,
\end{equation}
for $0 < s < t$. 

Now, we fix $s > 0$, $T > 0$, and we assume that $\omega_i(s) =   
\chi_i \nabla f \equiv \nabla(\chi_i f) - (\nabla \chi_i)f$ 
for some $f \in L^1(\real^2)$. 
Using \reff{SNint} together with the bounds \reff{estimSN1},
\reff{estimSN2} (for one vortex), we obtain for $s < t < s + T$:
\begin{eqnarray*}
  \|\omega_i(t)\|_{L^1} &\le& \frac{K_7}{(t{-}s)^{\frac12}} 
    \Bigl(\frac{t}{s}\Bigr)^\gamma \|\chi_i f\|_{L^1} 
    + K_6 \|(\nabla \chi_i)f\|_{L^1} 
    + \int_s^t K_6 \|Q_i(t')\omega(t')\|_{L^1} \d t' \\ 
  &+& \int_s^t \frac{K_7}{(t{-}t')^{\frac12}} \Bigl(\frac{t}
    {t'}\Bigr)^\gamma \|R_i(t')\omega(t')\|_{L^1} \d t'~.
\end{eqnarray*}
Summing over $i \in \{0,\dots,N\}$, we thus find
\begin{eqnarray*}
  \|\omega(t)\|_{L^1} &\le& \frac{K_7}{(t{-}s)^{\frac12}} 
    \Bigl(\frac{t}{s}\Bigr)^\gamma \|f\|_{L^1} 
    + \frac{C_1 K_6}{d}\|f\|_{L^1}  
    + \int_s^t K_6 \|Q(t')\|_{L^\infty}\|\omega(t')\|_{L^1} \d t'\\
  &+& \int_s^t \frac{K_7}{(t{-}t')^{\frac12}} \Bigl(\frac{t}
    {t'}\Bigr)^\gamma \|R(t')\|_{L^\infty}\|\omega(t')\|_{L^1} \d t'~.
\end{eqnarray*}
Thus, if we define (for a fixed $s > 0$)
$$
  \|\omega\|_T \,=\, \sup_{s < t < s+T} (t-s)^{\frac12} 
  \Bigl(\frac{s}{t}\Bigr)^\gamma \|\omega(t)\|_{L^1}~,
$$
we obtain using \reff{RQbdd}
$$
   \|\omega\|_T \,\le\, C_3\Bigl(1+\frac{T^{\frac12}}{d}\Bigr)
   \|f\|_{L^1} + C_4 \frac{T}{d^2} \|\omega\|_T + 
   C_5 \frac{T^{\frac12}}{d} \|\omega\|_T~,
$$
where $C_3, C_4, C_5 > 0$ are independent of $d$. If we now assume 
that $T$ is sufficiently small so that
$$
   \frac{T}{d^2} \,\le\, \min\Bigl\{1\,,\,\frac1{4C_4}\,,\,
   \frac1{(4C_5)^2}\Bigr\}~,
$$
then $\|\omega\|_T \le 4 C_3 \|f\|_{L^1}$, hence
$$
   \|S_N(t,s)\nabla f\|_{L^1} \,\le\, \frac{4C_3}{(t{-}s)^{\frac12}}
   \Bigl(\frac{t}{s}\Bigr)^\gamma \|f\|_{L^1}~, \quad
   s < t < s + T~.
$$
This concludes the proof of Proposition~\ref{estimateSN}. \QED

\subsection{Proof of Proposition~\ref{converge}}\label{proof63}

The proof follows the approach of \cite{gallay-wayne3}. Using
parabolic regularization, we first show that the trajectory 
$\{w_i(\tau)\}$ is relatively compact in $L^2(m)$. We next 
prove that the $\alpha$-limit set $\cA_i$ of this trajectory 
is fully invariant under the evolution defined by the autonomous
equation
\begin{equation}\label{wieqlimit}
  \frac{\partial w_i}{\partial \tau}(\xi,\tau) + 
  v_i(\xi,\tau)\cdot\nabla w_i(\xi,\tau) \,=\, (\cL w_i)(\xi,\tau)~,
\end{equation}
which is obtained by setting $R_i = 0$ in \reff{wieq}. Using
the main result of \cite{gallay-wayne3} we conclude that 
$\cA_i = \{\alpha_i G\}$, which proves the claim. 
We start with the compactness result: 

\begin{lemma}\label{compact}
For any $i \in \{1,\dots,N\}$ and any $m > 1$, the 
trajectory $\{w_i(\tau)\}_{\tau < \log T}$ is relatively 
compact in $ L^2(m)$.
\end{lemma}

\proof Fix $i \in \{1,\dots,N\}$. By \reff{wibdd}, for any $m > 1$, 
there exists $C_m > 0$ such that $\|w_i(\tau)\|_{L^2(m)} \le C_m$ 
for all $\tau < \log(T)$. Let $H^1(m)$ be the weighted Sobolev space
defined by
\begin{equation}\label{H1mdef}
  H^1(m) \,=\, \{w \in L^2(m) \,|\, \nabla w \in L^2(m)\}~, \quad
  \|w\|_{H^1(m)}^2 \,=\, \|w\|_{L^2(m)}^2 + \|\nabla w\|_{L^2(m)}^2~. 
\end{equation}
By Rellich's criterion, the inclusion $H^1(m{+}1) \hookrightarrow 
L^2(m)$ is compact. Therefore, to prove Lemma~\ref{compact}, 
it is sufficient to verify that $\nabla w_i(\tau)$ is bounded in 
$L^2(m)$ for all $m > 1$. To this end, we proceed as in 
(\cite{gallay-wayne3}, Lemma~2.1). Fix $m > 1$, $\tau_0 < \log(T)$, 
and consider the integral equation
$$
  \nabla w_i(\tau) \,=\,  \nabla S(\tau - \tau_0) w_i (\tau_0) 
  - \int_{\tau_0}^\tau \nabla S(\tau - s)(v_i(s)\cdot\nabla w_i(s) 
  + R_i(s)\cdot \nabla w_i(s)) \d s~,
$$
where $R_i$ is given by \reff{Rdef}. If $0 < \tau - \tau_0 \le 1$, 
we can bound, using Proposition~\ref{semigroup},
$$
  \|\nabla w_i(\tau)\|_{L^2(m)} \,\le\, \frac{K_3}
  {a(\tau{-}\tau_0)^{\frac12}}\|w_i(\tau_0)\|_{L^2(m)}
  + \int_{\tau_0}^\tau \frac{K_4}{a(\tau{-}s)^{\frac1p}}
  \|(v_i(s) + R_i(s))\cdot \nabla w_i(s)\|_{L^p(m)}
  \d s~,
$$
where $1 < p < 2$. If $\frac{1}{q} = \frac{1}{p} - \frac12$, we have 
according to~\reff{BS2}
$$
  \|v_j(s)\|_{L^q} \,\le\, C \|w_j(s)\|_{L^p} \,\le\,  
  C \|w_j(s)\|_{L^2(m)} \,\le\, C C_m~,
$$
and according to \reff{tdbdd}
$$
  \|\e^{\frac{\tau}{2}}\tilde u_0(\xi \e^{\frac{\tau}{2}}+z_i,
  \e^\tau)\|_{L^q_\xi}
  \,=\, \e^{\tau(\frac12-\frac1q)} \|\tilde u_0(\cdot,\e^\tau)\|_{L^q}
  \,=\, t^{\frac12-\frac1q} \|\tilde u_0(\cdot,t)\|_{L^q} 
  \,\le\, C_1 \|\mu_0\|~.
$$
Therefore, we find
$$
  \|(v_i(s) + R_i(s))\cdot \nabla w_i(s)\|_{L^p(m)} 
  \,\le\, \|v_i(s) + R_i(s)\|_{L^q} \|\nabla w_i(s)\|_{L^2(m)}
  \,\le\, C \|\nabla w_i(s)\|_{L^2(m)}~,
$$
hence
\begin{equation}\label{temp}
  \|\nabla w_i(\tau)\|_{L^2(m)} \,\le\, \frac{C_2}
  {a(\tau{-}\tau_0)^{\frac12}} + \int_{\tau_0}^\tau 
  \frac{C_3}{a(\tau{-}s)^{\frac1p}}\|\nabla w_i(s)\|_{L^2(m)}\d s~,
\end{equation}
where $C_2, C_3 > 0$ are independent of $\tau_0$. Now, choose 
$\tilde T > 0$ small enough so that
$$
   \sup_{\tau_0 < \tau \le \tau_0{+}\tilde T} 
   \int_{\tau_0}^\tau \frac{C_3 a(\tau{-}\tau_0)^{\frac12}}
   {a(\tau{-}s)^{\frac1p}a(s{-}\tau_0)^{\frac12}}\d s 
   \,\equiv\, \sup_{0 < \tau \le \tilde T} 
   \int_{0}^\tau \frac{C_3 a(\tau)^{\frac12}}
   {a(\tau{-}s)^{\frac1p}a(s)^{\frac12}}\d s 
   \,\le\, \frac12~\cdotp
$$
Then \reff{temp} implies that $\|\nabla w_i(\tau)\|_{L^2(m)} \le 
2C_2 a(\tau{-}\tau_0)^{-\frac12}$ for $\tau_0 < \tau < \min(
\tau_0{+}\tilde T,\log(T))$. Since $\tau_0 < \log(T)$ 
was arbitrary and since~$\tilde T$ is independent of~$\tau_0$, 
there exists~$C_4 > 0$ such that~$\|\nabla w_i(\tau)\|_{L^2(m)} \le 
C_4$ for all~$\tau < \log(T)$. \QED

\medskip
We next show that the term $R_i(\xi,\tau)$ in \reff{wieq} 
is negligible as $\tau \to -\infty$:

\begin{lemma}\label{lemremainder}
For any $i \in \{1,\dots,N\}$, any $m > 1$, and any $p \in (1,2)$,
the following holds:
$$
  \lim_{\tau \rightarrow - \infty}\|R_i(\tau) w_i(\tau)\|_{L^p(m)} = 0~.
$$
\end{lemma}

\proof Let $q \in (2,\infty)$ be such that $\frac{1}{q} = \frac{1}{p}
 - \frac12$. If $b(\xi) = (1{+}|\xi|^2)^{\frac12}$, we have for all 
$\tau < \log(T)$
$$
   \|R_i (\tau) w_i (\tau)\|_{L^p(m)} \,=\, \|b^m R_i (\tau) 
   w_i (\tau)\|_{L^p} \,\le\, \|b^{-1}R_i(\tau)\|_{L^q} 
   \|b^{m+1}w_i(\tau)\|_{L^2} \,\le\, C_1 \|b^{-1}R_i(\tau)\|_{L^q}~,
$$
where $C_1 > 0$ is independent of $\tau$. We claim that the last term 
in the right-hand side converges to zero as $\tau \to -\infty$. 
Indeed, if $0 < \nu < 1 - \frac2q$ and $j \neq i$, we have
$$
  \|b^{-1} v_j(\cdot - (z_j{-}z_i)\e^{-\frac{\tau}{2}},\tau)\|_{L^q}
  \,\le\, C \|b^{-1} b(\cdot - (z_j{-}z_i)\e^{-\frac{\tau}{2}})^{-\nu}
  \|_{L^\infty} \|b^\nu v_j(\tau)\|_{L^q}~.
$$
The first factor in the right-hand side is $\cO(\e^{\nu\frac{\tau}{2}})$ 
as $\tau \to -\infty$, and using \reff{weightedHLS} the second one
can be bounded by $C \|w_j(\tau)\|_{L^2(m)} \le C_2$, where $C_2 > 0$
is independent of $\tau$. On the other hand, applying 
Lemma~\ref{estimatetildeomega0} with $\chi(r) = 
(1{+}r)^{-\frac12}$, we obtain (with $t = \e^\tau$)
$$
  \|b^{-1}\e^{\frac{\tau}{2}}\tilde u_0(\xi \e^{\frac{\tau}{2}}+z_i,
  \e^\tau)\|_{L^q_\xi}
  \,=\, t^{\frac12-\frac1q} \Bigl\|\tilde u_0(x,t) 
  \Bigl(1 + \frac{|x-z_i|^2}{t}\Bigr)^{-\frac12}\Bigr\|_{L^q_x} 
  \,\build\hbox to 10mm{\rightarrowfill}_{t \to 0}^{}\, 0~,
$$
Thus $\|b^{-1}R_i(\tau)\|_{L^q} \to 0$ as $\tau \to -\infty$. \QED

\medskip
Now, fix $m > 1$ and let $\cA_i$ be the $\alpha$-limit set
in $L^2(m)$ of the trajectory $\{w_i(\tau)\}_{\tau < \log(T)}$. 
As is well known, $\cA_i$ is nonempty, compact, and attracts 
$w_i(\tau)$ in the sense that~$\dist_{L^2(m)} (w_i(\tau), \cA_i) 
\to 0$ as~$\tau \to -\infty$. Let also ${\Phi(\tau)}_{\tau\ge0}$ be 
the semiflow in $L^2(m)$ defined by the limiting equation 
\reff{wieqlimit}, see (\cite{gallay-wayne1}, Theorem~3.2).
Our last result is: 

\begin{lemma}\label{fullyinvariant}
For any $i \in \{1,\dots,N\}$ and any $\tau \ge 0$, we have
$\Phi (\tau) \cA_i =  \cA_i$.
\end{lemma}

\proof
It is clearly enough to prove the result for $0 \le \tau \le 1$. 
If $w_\infty \in \cA_i$, there exists a sequence~$\tau_n$ going
to~$-\infty$ such that~$\|w_i(\tau_n)-w_\infty\|_{L^2(m)} \to 0$ 
as~$n \to \infty$. From \reff{wiinteq2} we have
$$
 w_i (\tau_n + \tau) \,=\, S(\tau) w_i(\tau_n) -\int_0^\tau
 \e^{-\frac12(\tau - \tau')} \nabla\cdot S(\tau -\tau') 
 (v_i w_i + R_i w_i)(\tau_n +  \tau') \d \tau'~, 
$$
for $\tau \in [0,1]$. On the other hand, $W_i(\tau) \eqdef 
\Phi(\tau) w_\infty$ satisfies
$$
 W_i (\tau) \,=\, S(\tau) w_\infty  - \int_0^\tau  
 \e^{-\frac12(\tau - \tau')}\nabla\cdot S(\tau - \tau')
 (V_i W_i)(\tau') \d \tau'~, \quad 0 \le \tau \le 1~,
$$
where $V_i(\tau)$ is the velocity field obtained from $W_i(\tau)$
via the Biot-Savart law. Now we compute the difference of both 
expressions. Using Proposition~\ref{semigroup} and proceeding as in 
(\cite{gallay-wayne1}, Lemma~3.1) we obtain, for any $p \in (1,2)$, 
\begin{eqnarray*}
  &&\|w_i(\tau_n + \tau) - W_i (\tau)\|_{L^2(m)} \,\le\, K_3  
    \|w_i(\tau_n) - w_\infty\|_{L^2(m)}  \\
  &&\qquad + K_4\int_0^\tau \e^{-\frac12 (\tau - \tau')} 
    \frac{1}{a(\tau-\tau')^{\frac{1}{p}}} \|R_i(\tau_n + \tau') 
    w_i(\tau_n + \tau')\|_{L^p(m)} \d \tau' \\
  &&\qquad + K_4 \int_0^\tau \e^{-\frac12 (\tau - \tau')} 
    \frac{C}{a(\tau - \tau')^{\frac{1}{p}}} \left(
    \|w_i(\tau_n + \tau')\|_{L^2(m)} + \|W_i( \tau')\|_{L^2(m)} 
    \right) \\
  && \hspace{5.6cm} \times\, \|w_i(\tau_n  + \tau') -  W_i(\tau')\|_{L^2(m)}
    \d \tau'~.
\end{eqnarray*}
The first term in the right-hand side converges to zero as 
$n \to \infty$, and so does the second one (uniformly in 
$\tau \in [0,1]$) by Lemma~\ref{lemremainder}. Thus using the
uniform bound on $\|w_i(\tau)\|_{L^2(m)}$ we deduce that
$$
  \|w_i(\tau_n + \tau) - W_i (\tau)\|_{L^2(m)} \,\le\,
  \epsilon(n) + C \int_0^\tau  \frac{1}{a(\tau-\tau')^{\frac{1}{p}}}
  \|w_i(\tau_n + \tau') - W_i(\tau')\|_{L^2(m)} \d \tau'~, 
$$
where $\epsilon(n) \to 0$ as $n \to \infty$ uniformly in $\tau 
\in [0,1]$, and where $C > 0$ is independent of $n$ and $\tau$. 
Gronwall's lemma then implies
$$
  \lim_{n \rightarrow \infty }\sup_{0 \le \tau \le 1} \|w_i(\tau_n
  + \tau) -  W_i(\tau)\|_{L^2(m)} \,=\, 0~.
$$
In particular $W_i(\tau) \equiv \Phi(\tau)w_\infty \in \cA_i$ 
for any $\tau \in [0,1]$, hence $\Phi(\tau) \cA_i \subset 
\cA_i$ for any $\tau\in [0,1]$. 

Conversely let $w_\infty \in \cA_i$ and $0 \le \tau \le 1$. 
There exists a sequence~$\tau_n$ going to~$-\infty$ such 
that~$\|w_i(\tau_n) - w_\infty\|_{L^2(m)} \to 0$ as~$n \to \infty$.
Up to the extracting a subsequence, we can suppose 
that~$w_i(\tau_n - \tau)$ converges in $L^2(m)$ towards some~$W_\infty 
\in \cA_i$ as~$n$ goes to infinity. By the argument above, we 
have~$\Phi(\tau) W_\infty = w_\infty$, hence~$\cA_i \subset \Phi (\tau)
\cA_i$. \QED

\medskip In (\cite{gallay-wayne3}, Lemma~3.3) it is shown that, 
if $\cA$ is a bounded subset of $L^2(m)$ satisfying 
$\Phi(\tau)\cA = \cA$ for all $\tau \ge 0$, then necessarily
$\cA \subset \{\alpha G\,|\,\alpha \in \real\}$. Applying this
result to the~$\alpha$-limit set of~$\{w_i(\tau)\}$, we 
obtain~$\cA_i = \{\alpha_i G\}$, since any~$w \in \cA_i$ 
satisfies~$\inttwo w(\xi)\d\xi = \alpha_i$. This concludes the proof 
of Proposition~\ref{converge}. \QED

\subsection{Proof of Proposition~\ref{talphai}}\label{proof64}

Estimates \reff{talpha1} and \reff{talpha2} are established
in (\cite{gallay-wayne3}, Section 4.2). To prove iii), 
we first observe that it is sufficient to establish 
\reff{talpha3} for $0 < \tau \le \tau_0$, where $\tau_0 > 0$ 
is arbitrary. Indeed, once this is done, we have for $\tau > \tau_0$:
\begin{eqnarray*}
  \|T_\alpha(\tau)\nabla w\|_{L^2(m)} &=&  
  \|T_\alpha(\tau{-}\tau_0)T_\alpha(\tau_0)\nabla w\|_{L^2(m)} \\
  &\le& K_8 \,\e^{-\frac{(\tau-\tau_0)}{2}}\|T_\alpha(\tau_0)
  \nabla w\|_{L^2(m)} \,\le\, \frac{K_8 K_9 \,\e^{-\frac{\tau}{2}}}
  {a(\tau_0)^\frac1q} \|w\|_{L^q(m)}~.
\end{eqnarray*}
To prove \reff{talpha3} for small $\tau$, we consider the 
auxiliary equation in $L^2(m)^2$:
\begin{equation}\label{auxf}
  \frac{\partial f}{\partial \tau} + \alpha \Bigl(v^G \div f + 
  v^{\div f} G\Bigr) \,=\, \Bigl(\cL - \frac12\Bigr)f~.
\end{equation}
Here, $f(x,t) \in \real^2$ is a vector field, and $v^{\div f}$
denotes the velocity field obtained from the scalar~$\div f$ 
via the Biot-Savart law~\reff{BS2}. As we shall see, this equation 
defines a semigroup in~$L^2(m)^2$, which we denote by~$\hat 
T_\alpha(\tau)$. A straightforward calculation shows 
that the semigroups~$T_\alpha(\tau)$ and $\hat T_\alpha(\tau)$ 
are related via
$$
  T_\alpha(\tau)\div f \,=\, \div (\hat T_\alpha(\tau)f)~, \quad
  f \in H^1(m)^2~, \quad \tau \ge 0~,
$$
where $H^1(m)$ is defined in \reff{H1mdef}. Assertion iii) in 
Proposition~\ref{talphai} is now a direct consequence of this
identity and of the following result: 

\begin{lemma}\label{auxT}
Equation~\reff{auxf} defines a strongly continuous semigroup
$\hat T_\alpha(\tau)$ in $L^2(m)^2$ for any~$m > 1$. If~$q \in [1,2]$
and~$\tau > 0$, $\hat T_\alpha(\tau)$ can be
extended to a bounded operator from~$L^q(m)^2$ to~$H^1(m)^2$, 
and there exist~$\tau_0 > 0$ and~$C > 0$ such that
\begin{equation}\label{talpha4}
  \|\hat T_\alpha(\tau)f\|_{L^2(m)} \,\le\, \frac{C}
    {a(\tau)^{\frac1q-\frac12}} \|f\|_{L^q(m)}~, \quad
  \|\nabla\hat T_\alpha(\tau)f\|_{L^2(m)} \,\le\, \frac{C}
    {a(\tau)^{\frac1q}} \|f\|_{L^q(m)}~,
\end{equation}
for $\tau \in (0,\tau_0]$. 
\end{lemma}

\noindent{\bf Proof.} We consider the integral equation associated
with \reff{auxf}, namely
\begin{equation}\label{intauxf}
  f(\tau) \,=\, \e^{-\frac12\tau}S(\tau)f_0 - \alpha
  \int_0^\tau \e^{-\frac12(\tau{-}s)}S(\tau{-}s)
  \Bigl(v^G \div f(s) + v^{\div f(s)}G\Bigr)\d s~.
\end{equation}
We assume that $f_0 \in L^q(m)^2$ for some $q \in [1,2]$ and 
some $m > 1$. Given $\tau_0 > 0$, we shall solve~\reff{intauxf}
in the Banach space $X \eqdef \{f \in C^0((0,\tau_0],H^1(m)^2)\,|\,
\|f\|_X < \infty\}$, where
$$
  \|f\|_X \,=\, \sup_{0 < \tau \le \tau_0} a(\tau)^{\frac1q-\frac12}
  \|f(\tau)\|_{L^2(m)} + \sup_{0 < \tau \le \tau_0} a(\tau)^{\frac1q}
  \|\nabla f(\tau)\|_{L^2(m)}~.
$$
Let $f \in X$ and denote by $F(\tau)$ the expression in the 
right-hand side of \reff{intauxf}. Using the estimates collected
in Proposition~\ref{semigroup}, we obtain for $\tau \in (0,\tau_0]$:
$$
   \|F(\tau)\|_{L^2(m)} \,\le\, \frac{K_6\,\e^{-\frac12\tau}}
   {a(\tau)^{\frac1q-\frac12}} \|f_0\|_{L^q(m)} + 
   |\alpha| \int_0^\tau \frac{K_6\,\e^{-\frac12(\tau{-}s)}}
   {a(\tau{-}s)^{\frac1p-\frac12}} \|v^G \div f(s) + 
   v^{\div f(s)}G\|_{L^p(m)}\d s~,
$$
where $1 < p < 2$. If $\frac{1}{p'} = \frac{1}{p} - \frac12$, we
obtain using \reff{HLS} and H\"older's inequality
\begin{eqnarray*}
  &&\|v^G \div f(s)\|_{L^p(m)} \,\le\, \|v^G\|_{L^{p'}}
   \|\div f(s)\|_{L^2(m)} \,\le\, C \|G\|_{L^p} \|\div f(s)\|_{L^2(m)}~,\\
  &&\|v^{\div f(s)} G\|_{L^p(m)} \,\le\, C \|v^{\div f(s)}\|_{L^{p'}}
   \|G\|_{L^2(m)} \,\le\, \|\div f(s)\|_{L^p} \|G\|_{L^2(m)}~,
\end{eqnarray*}
hence both terms can be bounded by $C\|G\|_{L^2(m)}\|\div f(s)\|_{L^2(m)}$.
Therefore, 
$$
  {a(\tau)^{\frac1q-\frac12}} \|F(\tau)\|_{L^2(m)} \,\le\, 
  C\e^{-\frac12\tau} \|f_0\|_{L^q(m)} + C |\alpha| \Lambda_1(\tau)
  \|f\|_X~,
$$
where
$$
   \Lambda_1(\tau) \,=\, \int_0^\tau \e^{-\frac12(\tau{-}s)}
   \frac{a(\tau)^{\frac1q-\frac12}}{a(\tau{-}s)^{\frac1p-\frac12}
   a(s)^{\frac1q}}\d s~.
$$
Using similar estimates, one finds
$$  
  {a(\tau)^{\frac1q}} \|\nabla F(\tau)\|_{L^2(m)} \,\le\, 
  C\e^{-\frac12\tau} \|f_0\|_{L^q(m)} + C |\alpha| \Lambda_2(\tau)
  \|f\|_X~,
$$
where
$$
   \Lambda_2(\tau) \,=\, \int_0^\tau \e^{-\frac12(\tau{-}s)}
   \frac{a(\tau)^{\frac1q}}{a(\tau{-}s)^{\frac1p}
   a(s)^{\frac1q}}\d s~.
$$
Thus, there exist positive constants $C_1$, $C_2$ such that, 
for all $f \in X$, 
$$
   \|F\|_X \,\le\, C_1 \|f_0\|_{L^q(m)} + C_2 |\alpha| 
   \Lambda(\tau_0) \|f\|_X~,
$$
where
$$
  \Lambda(\tau_0) \,=\, \sup_{0 < \tau \le \tau_0} \Lambda_1(\tau)
  + \sup_{0 < \tau \le \tau_0} \Lambda_2(\tau)~.
$$
If we now choose $\tau_0 > 0$ small enough so that $C_2 |\alpha| 
\Lambda(\tau_0) \le \frac12$, it follows from these estimates that
\reff{intauxf} has a unique solution $f \in X$ such that 
$\|f\|_X \le 2 C_1 \|f_0\|_{L^q(m)}$. This proves \reff{talpha4}.
\QED

%%%%%%%%%%%%%%%%%%%%%%%%%%%%%%%%%%%%%%%%%%%%%%%%%%%%%%%%%%%%%%%%%%%%%%%

\bibliographystyle{plain}

\begin{thebibliography}{10}

\bibitem{ben-artzi}
M.~Ben-Artzi.
Global solutions of two-dimensional Navier-Stokes and Euler
equations.
{\em Arch. Rational Mech. Anal.} {\bf 128} (1994), 329--358.

\bibitem{brezis}
H.~Brezis.
Remarks on the preceding paper by M. Ben-{A}rtzi: ``Global
solutions of two-dimensional Navier-Stokes and Euler
equations''.
{\em Arch. Rational Mech. Anal.} {\bf 128} (1994), 359--360.

\bibitem{cannone-book}
M.~Cannone. 
{\sl Ondelettes, paraproduits et Navier-Stokes.}  
Diderot Editeur, Paris, 1995. 

\bibitem{cannone-planchon}
M.~Cannone and F.~Planchon.
Self-similar solutions for Navier-Stokes equations in $\real^3$.
{\em Comm. Partial Differ. Equations} {\bf 21} (1996), 179--193.

\bibitem{carlen-loss}
E.~A. Carlen and M.~Loss.
Optimal smoothing and decay estimates for viscously damped
conservation laws, with applications to the {$2$}-D Navier-Stokes
equation.
{\em Duke Math. J.} {\bf 81} (1995), 135--157 (1996). 

\bibitem{carpio}
A.~Carpio.
Asymptotic behavior for the vorticity equations in dimensions two and
three.
{\em Comm. Partial Differential Equations} {\bf 19} (1994), 827--872.

\bibitem{chorin-marsden}
A.~Chorin and J.~Marsden.  
{\sl A mathematical introduction to fluid mechanics.}
Second edition. Texts in Applied Mathematics, {\bf 4}. 
Springer-Verlag, New York, 1990. 

\bibitem{cottet}
G.-H. Cottet.
Equations de Navier-Stokes dans le plan avec tourbillon initial
mesure. 
{\em C. R. Acad. Sci. Paris S\'er. I Math.} {\bf 303} (1986), 105--108.

\bibitem{gp} I. Gallagher and F. Planchon. 
On global infinite energy solutions to the Navier--Stokes equations 
in two dimensions.
{\em Arch. Rational Mech. Anal.} {\bf 161} (2002),
307--337. 

\bibitem{gallay-wayne1}
Th. Gallay and C.~E. Wayne.
Invariant manifolds and the long-time asymptotics of the
Navier-Stokes and vorticity equations on {$\real^2$}.
{\em Arch. Rational Mech. Anal.} {\bf 163} (2002), 209--258.

\bibitem{gallay-wayne2}
Th. Gallay and C.~E. Wayne.
Long-time asymptotics of the Navier-Stokes and vorticity
equations on {$\real^3$}.
{\em R. Soc. Lond. Philos. Trans. Ser. A Math. Phys. Eng. Sci.}
{\bf 360} (2002), 2155--2188. 

\bibitem{gallay-wayne3}
Th. Gallay and C.~E. Wayne.
Global stability of vortex solutions of the two-dimensional
Navier-Stokes equation. {\em Comm. Math. Phys.}, to appear.
Preprint version available at {\tt http://www.arXiv.org/math.AP/0402449}.

\bibitem{germain} P. Germain. 
Existence globale de solutions de l'\'equation de Navier-Stokes 
2D avec donn\'ees initiales dans $\partial\bmo$.
In preparation.

\bibitem{giga-miyakawa1}
Y.~Giga and T.~Miyakawa. 
Solutions in $L^r$ of the Navier-Stokes initial value problem. 
{\em Arch. Rational Mech. Anal.} {\bf 89} (1985), 267--281.

\bibitem{giga-kambe}
Y.~Giga and T.~Kambe.
Large time behavior of the vorticity of two-dimensional viscous
flow and its application to vortex formation.
{\em Comm. Math. Phys.} {\bf 117} (1988), 549--568.

\bibitem{GMO}
Y.~Giga, T.~Miyakawa, and H.~Osada.
Two-dimensional Navier-Stokes flow with measures as initial
vorticity.
{\em Arch. Rational Mech. Anal.} {\bf 104} (1988), 223--250.

\bibitem{giga-miyakawa2}
Y.~Giga and T.~Miyakawa.
Navier-Stokes flow in $\real^3$ with measures as initial vorticity 
and Morrey spaces. 
{\em Comm. Partial Differential Equations} {\bf 14} (1989), 577--618.

\bibitem{gigabook}
M.-H.~Giga and Y.~Giga.
{\sl Nonlinear partial differential equations -- asymptotic 
behaviour of solutions and self-similar solutions}.
Book in preparation. 

\bibitem{henry}
D.~Henry.
{\em Geometric theory of semilinear parabolic equations}.
Springer-Verlag, Berlin, 1981.

\bibitem{kato84}
T.~Kato.
Strong $L^p$-solutions of the Navier-Stokes equation in $\real^m$, 
with applications to weak solutions. 
{\em Math. Z.} {\bf 187} (1984), 471--480.

\bibitem{kato}
T.~Kato.
The Navier-Stokes equation for an incompressible fluid in
$\real^2$ with a measure as the initial vorticity.
{\em Differential Integral Equations} {\bf 7} (1994), 949--966.

\bibitem{kochtataru} H. Koch and D. Tataru.  
Well-posedness for the Navier--Stokes equations.
{\em Advances in Mathematics} {\bf 157} (2001), 22--35. 

\bibitem{ladyzhenskaia}
O.~Ladyzhenskaya,
{\sl The mathematical theory of viscous incompressible flow.} 
Second English edition. Mathematics and its Applications, 
Vol. 2. Gordon and Breach, New York-London-Paris, 1969.

\bibitem{leray} J.~Leray.
Etude de diverses \'equations int\'egrales non lin\'eaires et
de quelques probl\`emes que pose l'hydrodynamique.
{\em J. Math. Pures. Appl.} {\bf 12} (1933), 1--82. 

\bibitem{meyer} Y. Meyer. 
{\sl Wavelets, Paraproducts and Navier--Stokes}.  
Current Developments in Mathematics, International Press, Cambridge, 
Massachussets, 1996.

\bibitem{osada}
H.~Osada.
Diffusion processes with generators of generalized divergence form.
{\em J. Math. Kyoto Univ.} {\bf 27} (1987), 597--619.

\bibitem{stein} E. Stein. {\sl Harmonic Analysis}. 
Princeton University Press, 1993.

\bibitem{temam}
R.~Temam. {\sl Navier-Stokes equations. Theory and numerical
analysis}. 
Reprint of the 1984 edition. AMS Chelsea Publishing, Providence, RI, 2001. 

\bibitem{weissler}
F. Weissler. 
The Navier-Stokes initial value problem in $L^p$. 
{\em Arch. Rational Mech. Anal.} {\bf 74} (1980), 219--230.

\end{thebibliography}

\end{document}